\documentclass[a4paper,11pt,reqno]{amsart}

\usepackage[T1]{fontenc}
\usepackage{lmodern}
\usepackage{a4wide}

\usepackage{amsmath,amssymb,mathtools}
\allowdisplaybreaks[4]
\usepackage{mathrsfs}

\usepackage{enumitem}
\usepackage{needspace}

\usepackage{xcolor}
\usepackage[colorlinks=true,citecolor=teal,linkcolor=blue,urlcolor=blue]{hyperref}
\usepackage{orcidlink}

\makeatletter
\def\@maketitle{\normalfont\normalsize
  \let\@makefnmark\relax  \let\@thefnmark\relax
  \ifx\@empty\@date\else \@footnotetext{\@setdate}\fi
  \ifx\@empty\@subjclass\else \@footnotetext{\@setsubjclass}\fi
  \ifx\@empty\@keywords\else \@footnotetext{\@setkeywords}\fi
  \ifx\@empty\thankses\else \@footnotetext{\def\par{\let\par\@par}\@setthanks}\fi
  \global\topskip42\p@
  \@settitle
  \ifx\@empty\authors \else \@setauthors \fi
  \@setabstract
  \normalsize
  \dimen@34\p@ \advance\dimen@-\baselineskip
  \vskip\dimen@\relax
}
\def\@citestyle{\m@th\upshape\mdseries\color{teal}}
\def\@biblabel#1{\textcolor{teal}{\textbf{[#1]}}}
\makeatother

\numberwithin{equation}{section}

\theoremstyle{plain}

\newtheorem{theorem}{Theorem}[section]
\newtheorem{lemma}[theorem]{Lemma}
\newtheorem{proposition}[theorem]{Proposition}

\newtheorem{conjecture}[theorem]{Conjecture}

\theoremstyle{remark}

\newtheorem{remark}[theorem]{Remark}

\newcommand{\R}{\mathbb{R}}
\newcommand{\Sph}{\mathbb{S}}

\newcommand{\Ric}{\operatorname{Ric}}
\newcommand{\Rm}{\operatorname{Rm}}
\newcommand{\divg}{\operatorname{div}}
\newcommand{\trg}{\operatorname{tr}}
\newcommand{\Id}{\operatorname{Id}}
\newcommand{\II}{\mathrm{II}}

\newcommand{\inner}[2]{\left\langle #1,#2\right\rangle}
\newcommand{\dmu}{\,\mathrm{d}\mu}

\newcommand{\papertitle}{Rigidity of compact quasi-Einstein manifolds with boundary and low-dimensional focal sets}
\newcommand{\papershorttitle}{Quasi-Einstein manifolds with low-dimensional focal sets}

\title[\papershorttitle]{\papertitle}

\author[M. A. R. M. Hor\'acio]{Matheus Andrade Ribeiro de Moura Hor\'acio\,\orcidlink{0009-0001-5658-7467}}
\thanks{ORCID iD: \url{https://orcid.org/0009-0001-5658-7467}.}

\address{
Departamento de Matem\'atica,
Universidade de Bras\'ilia,
70910-900, Bras\'ilia, DF,
Brazil
}

\email{matheus.andrade5488@gmail.com}

\date{}



\hypersetup{
  pdftitle={\papertitle},
  pdfauthor={Matheus Andrade Ribeiro de Moura Horacio},
  pdfsubject={Differential geometry},
  pdfkeywords={Quasi-Einstein manifolds, rigidity, constant scalar curvature, focal set}
}

\begin{document}
\raggedbottom

\begin{abstract}
In this article, we investigate nontrivial compact $m$-quasi-Einstein manifolds with boundary, $m>1$ and constant scalar curvature. We first prove that such a manifold is rigid whenever the set of critical points of its potential function is totally geodesic, which extends to quasi-Einstein manifolds a theorem of Li, Li, Ou, Qu, Wu and Zhou for shrinking Ricci solitons. Then, we prove that, in every dimension, this set is totally geodesic whenever it has dimension two or three, so that the manifold is rigid in this case. The proof is based on a Bochner identity for the covariant derivative of the Ricci tensor at the points of this set. In particular, every such manifold of dimension five is rigid. Finally, we obtain rigidity in arbitrary dimension under an integral Weyl pinching condition when the set of critical points has dimension four.
\end{abstract}

\maketitle

\section{Introduction and main results}\label{sec_intro}

Let $(M^n,g)$, $n\geq2$, be a complete Riemannian manifold, possibly with boundary $\partial M$. We say that $(M^n,g,u,\lambda)$ is an $m$-quasi-Einstein manifold if $u$ is a smooth function on $M$ and $\lambda$ and $m>0$ are constants such that
\begin{align}\label{fundeq}
	\nabla^2u=\frac{u}{m}\big(\Ric-\lambda g\big),\qquad u>0\ \text{ in the interior of } M,\qquad u=0\ \text{ on }\partial M,
\end{align}
see \cite{csw,hpw1,hpw2}. It is well known that, when $m$ is an integer, these manifolds are the bases of the warped product Einstein metrics whose fiber has dimension $m$ \cite[Proposition 1.1]{hpw1}; see also \cite{besse}. We say that $(M^n,g,u,\lambda)$ is \emph{trivial} when $u$ is constant. In this case, $\partial M=\emptyset$ and $g$ is an Einstein metric. Following He, Petersen and Wylie \cite[Definition 1.2]{hpw2}, we say that it is \emph{rigid} when it is Einstein or its universal cover is a product of Einstein manifolds. Consider a nontrivial rigid $m$-quasi-Einstein manifold with boundary, $\lambda>0$ and $m>1$. By \cite[Proposition 2.6]{hpw2}, its universal cover is a product $\bar N^k\times\Sph^{q}_+(\alpha)$, where $\bar N^k$ is an Einstein manifold with $\Ric_{\bar N}=\lambda g_{\bar N}$ and $\Sph^q_+(\alpha)$ is a hemisphere of the round sphere of sectional curvature $\alpha=\lambda/(m+q-1)$. Moreover, $u$ is a multiple of the cosine of $\sqrt{\alpha}$ times the distance to the pole of the hemisphere. Here either factor may have dimension zero, and $\Sph^1_+(\alpha)$ denotes an interval of length $\pi/\sqrt{\alpha}$.

When $m$ is finite, equation \eqref{fundeq} is the counterpart of the gradient Ricci soliton equation
\begin{align}\label{eq_soliton}
	\Ric+\nabla^2f=\lambda g,
\end{align}
and, formally, we recover \eqref{eq_soliton} from \eqref{fundeq} if we set $u=e^{-f/m}$ and let $m\to\infty$; see \cite{csw}. In either case, the Ricci tensor determines the Hessian of the potential function. Let us look at the points where the gradient of the potential function does not vanish. There the gradient gives a privileged direction, and the level sets are hypersurfaces whose second fundamental form is, by the equation, a multiple of the restriction of $\lambda g-\Ric$. This direction is special, in the sense that hypotheses on the curvature in this direction alone may determine the whole manifold, as several structure and classification results show. For instance, it was proved by Petersen and Wylie \cite[Theorem 1.2]{pewy} that a complete gradient shrinking Ricci soliton is rigid if and only if its scalar curvature is constant and it is radially flat, that is, $\Rm(\nabla f,X,\nabla f,X)=0$ for every $X$. In other words, rigidity is decided by the sectional curvatures of the planes which contain this single direction. The corresponding rigidity theorem for quasi-Einstein manifolds is due to He, Petersen and Wylie \cite[Theorem 1.11]{hpw2}. Moreover, this direction splits off under conditions on the Weyl tensor, or on the level hypersurfaces. More precisely, around the regular points of the potential function, the metric is a warped product, or a multiply warped product, with one-dimensional base tangent to the gradient of the potential function and fibers tangent to its level hypersurfaces, whence the potential function depends only on the base. For quasi-Einstein manifolds with harmonic Weyl tensor $W$ and $W(\nabla u,\cdot,\cdot,\nabla u)=0$, this is \cite[Theorem 1.2]{hpw1}. The same was proved in \cite{bhs} for gradient almost Ricci solitons with harmonic Weyl tensor, and in \cite{horlevel} for gradient Ricci solitons whose level hypersurfaces have parallel second fundamental form.

The set where this privileged direction ceases to exist, namely the critical set of the potential function, is no less important, and its study is the complementary approach adopted in this article. We begin with shrinking solitons. Cao and Zhou \cite[Theorem 1.1]{caozhou} proved that the potential function $f$ of a complete noncompact gradient shrinking Ricci soliton grows quadratically with the distance to a fixed point. Thus $f$ is proper, and it attains its minimum. Suppose, in addition, that the scalar curvature $R$ is constant. Then the identity $R+|\nabla f|^2-2\lambda f=\mathrm{const}$ \cite{pewy} shows that $|\nabla f|^2$ is an affine function of $f$, and therefore the critical points of $f$ are exactly its minimum points. Moreover, $\Delta f=n\lambda-R$ is constant. Consequently, both $|\nabla f|^2$ and $\Delta f$ are functions of $f$, which means that $f$ is an isoparametric function in the sense of Wang \cite{wang}; see also \cite{tranyamazaki} for Ricci solitons with isoparametric potential functions. The theory of these functions applies, and we have the following. The minimum set of $f$ is a compact submanifold, called the focal set of $f$, and its dimension is $R/\lambda$ \cite{chengzhou,fr16}. The regular level sets of $f$ are the tubes around the focal set, and they have constant mean curvature. Furthermore, the focal set is a minimal submanifold \cite{getang,wang}. A rigid soliton is covered by the product of an Einstein manifold $\bar N$ and a Euclidean space, and its potential function is, up to a constant, $\frac{\lambda}{2}|x|^2$ on the Euclidean factor \cite{pewy}. Hence its focal set is the image of $\bar N\times\{0\}$ under the covering map, and it is totally geodesic. Li, Li, Ou, Qu, Wu and Zhou \cite[Corollary 1.1]{lloqwz} proved the converse, namely that a complete gradient shrinking Ricci soliton with constant scalar curvature is rigid whenever its focal set is totally geodesic.

Cao conjectured that every complete gradient shrinking Ricci soliton with constant scalar curvature is rigid in the sense of Petersen and Wylie \cite{pewy}; see \cite{chengzhou}. Rigidity of shrinking solitons is also known under several other curvature conditions; see, for instance, \cite{flgr11,ms13,pw10}. Concerning the conjecture, Fern\'andez-L\'opez and Garc\'ia-R\'io proved in \cite[Theorem 1]{fr16} that the scalar curvature of such a soliton belongs to $\{0,\lambda,\ldots,n\lambda\}$. They also proved that $R=\lambda$ does not occur and that $R=(n-1)\lambda$ implies rigidity \cite[Theorem 10]{fr16}. The values $R=0$ and $R=n\lambda$ correspond to the Gaussian soliton and to compact Einstein manifolds, respectively \cite[Propositions 3.1 and 3.3]{pewy}. In dimension four, the conjecture was proved by Cheng and Zhou \cite{chengzhou}, and Ou, Qu and Wu \cite{ouquwu} gave another proof of the case $R=2\lambda$; see also \cite{wwres}. In dimension five, the remaining values are $R=3\lambda$ and $R=2\lambda$. They were treated by Li, Ou, Qu and Wu in the preprints \cite{loqw1} and \cite{loqw2}, respectively, and by the author in \cite{horfocal}. For K\"ahler solitons, the conjecture holds in every dimension \cite{hloqwz,horkahler,lowz}, and we refer to \cite{lqw,tranyamazaki} for related results. In arbitrary dimension, the theorem of \cite{lloqwz} mentioned above reduces the conjecture to proving that the focal set is totally geodesic. It was proved in \cite[Theorem 1.1]{horfocal} that this is the case whenever the focal set has dimension two or three. The proof uses a Bochner identity for the covariant derivative of the Ricci tensor at the points of the focal set. Consequently, the conjecture of Cao holds when $R\leq3\lambda$ \cite[Corollary 1.3]{horfocal}.

The same set appears for quasi-Einstein manifolds, expressed through a different potential function. To see this, let $(M^n,g,u,\lambda)$ be a nontrivial compact $m$-quasi-Einstein manifold with boundary, and consider the function $f=-m\log u$, which is defined in the interior of $M$ and satisfies
\[
	\Ric+\nabla^2f-\frac1m\,\mathrm{d}f\otimes\mathrm{d}f=\lambda g,
\]
an equation which formally tends to \eqref{eq_soliton} as $m\to\infty$. The minimum set of $f$ is the set $\operatorname{MAX}(u)$ of the maximum points of $u$. The function $f$ tends to infinity at $\partial M$, as the potential function of a complete noncompact shrinking soliton does at infinity, while $u$ is smooth on the compact manifold $M$ and vanishes on $\partial M$, so that it attains its maximum $U$ in the interior of $M$. Now suppose that $m>1$ and that the scalar curvature is constant. Then
\[
	|\nabla u|^2=\alpha\big(U^2-u^2\big)
\]
for some constant $\alpha>0$; see \eqref{eq_transnormal}. Therefore, exactly as for shrinking solitons, the critical points of $u$ are its maximum points. On the other hand, $\Delta u=\frac um(R-n\lambda)$ is also a function of $u$. Thus the potential function $u$ is isoparametric \cite{crz}, and its focal set $F=\operatorname{MAX}(u)$ is a closed connected minimal submanifold contained in the interior of $M$. We also have $u=U\cos(\sqrt{\alpha}\,r)$, where $r$ is the distance to $F$. In particular, the level sets of $u$ are the sets of the points at a fixed distance from $F$; see Proposition \ref{prop_structure} and Lemma \ref{lem_first_order}.

For quasi-Einstein manifolds, the natural setting for the problem corresponding to the conjecture of Cao is that of compact manifolds with boundary, $\lambda>0$ and $m>1$. Indeed, a nontrivial complete $m$-quasi-Einstein manifold is compact if and only if $\lambda>0$ \cite[Theorem 4.1]{hpw1}. Moreover, by \cite[Proposition 2.1]{csw}, a compact $m$-quasi-Einstein manifold without boundary and with constant scalar curvature is trivial. As for the case $m=1$, the scalar curvature is always constant under the usual additional assumption $\Delta u=-\lambda u$, see \cite[p. 160]{hpw2}, and therefore the hypothesis of constant scalar curvature gives no information. Finally, for $\lambda<0$, He, Petersen and Wylie constructed complete noncompact four-dimensional examples with constant scalar curvature which are not rigid \cite[Theorem 1.5]{hpw2}.

Let $(M^n,g,u,\lambda)$ be a nontrivial compact $m$-quasi-Einstein manifold with boundary, $m>1$ and constant scalar curvature $R$. It was proved by Costa, Ribeiro Jr and Zhou \cite[Theorem 1 and Proposition 5]{crz} that
\begin{align}\label{eq_quantization}
	R=\frac{k(m-n)+n(n-1)}{m+n-k-1}\,\lambda\qquad\text{for some } k\in\{0,2,3,\ldots,n-1\},
\end{align}
where $k$ is the dimension of $F$. This result is the counterpart of the theorem of Fern\'andez-L\'opez and Garc\'ia-R\'io mentioned above. Moreover, the right-hand side of \eqref{eq_quantization} tends to $k\lambda$ as $m\to\infty$. The following cases are known. If $k=0$, then $M$ is Einstein, and it is isometric to a hemisphere \cite[Proposition 2.4]{hpw2}. If $k=n-1$, then the universal cover of $M$ is the product of an interval and a compact Einstein manifold \cite[Theorem 2]{crz}. In dimension three, these are the only possible values of $k$, and rigidity was proved by He, Petersen and Wylie \cite[Theorem 1.3]{hpw2}; see also \cite[Corollary 1]{crz}. In dimension four, the only remaining value is $k=2$, which was settled by Costa, Ribeiro Jr and Zhou \cite[Theorem 3]{crz} with a proof inspired by \cite{chengzhou}. In dimension five, the remaining values are $k=2$ and $k=3$. The case $k=3$ was considered by Z. Cao in the preprint \cite[Theorem 1.6]{zcao}, whose proof is inspired by \cite{loqw1}. We refer to \cite{cls,crs,dgr} for related results on quasi-Einstein manifolds, and to \cite{crs} for the noncompact case. Very recently, Baltazar, Carvalho and Chen \cite{bcc} proved that a compact quasi-Einstein manifold with boundary of dimension $n\geq5$ and zero radial Weyl curvature is isometric, up to scaling, to the standard hemisphere or to a warped product $\mathrm{d}t^2+\psi^2g_L$, where $g_L$ is an Einstein metric with nonnegative Ricci curvature. It is worth mentioning that this hypothesis is not satisfied by the rigid manifolds $\bar N^k\times\Sph^q_+(\alpha)$ with $k\geq2$ and $q\geq2$. Indeed, if $X$ is a unit vector tangent to $\bar N$ and $\nu$ is a unit vector tangent to the hemisphere and parallel to $\nabla u$, then the Weyl tensor $W$ of the product satisfies
\[
	W(X,\nu,X,\nu)=-\frac{(q-1)\big(\lambda+(k-1)\alpha\big)}{(n-1)(n-2)}\neq0.
\]
In view of the results above, it is natural to formulate the following counterpart of the conjecture of Cao.

\begin{conjecture}\label{conj_main}
Every nontrivial compact $m$-quasi-Einstein manifold with boundary, $m>1$ and constant scalar curvature is rigid.
\end{conjecture}

To the best of our knowledge, He, Petersen and Wylie \cite[p. 160]{hpw2} were the first to ask whether quasi-Einstein manifolds with constant scalar curvature can be classified. For compact manifolds with boundary, Costa, Ribeiro Jr and Zhou asked in the introduction of \cite{crz} whether such a manifold must be one of the known examples; see also \cite{crs,zcao}. By the results mentioned above, Conjecture \ref{conj_main} holds when $n\leq4$. Thus the open cases are $2\leq k\leq n-2$ in dimension $n\geq5$.

The main goal of this article is to investigate Conjecture \ref{conj_main} by means of the focal set $F$. To this end, we extend to quasi-Einstein manifolds the method of \cite{horfocal}. In what follows, $q=n-k$ denotes the codimension of $F$ and $\II$ its second fundamental form. With this notation, $\alpha=\lambda/(m+q-1)$ by Proposition \ref{prop_structure}. On a rigid manifold, $F$ is the image, under the universal covering map, of the product of $\bar N^k$ and the pole of the hemisphere. In particular, it is totally geodesic; see also \cite[Proposition 3.14]{hpw2}. Thus the vanishing of $\II$ is a necessary condition for rigidity, which only involves a compact submanifold contained in the interior of $M$. Our first result shows that this condition is also sufficient.

\begin{theorem}\label{thm_criterion}
Let $(M^n,g,u,\lambda)$ be a nontrivial compact $m$-quasi-Einstein manifold with boundary, $m>1$ and constant scalar curvature. If $\operatorname{MAX}(u)$ is totally geodesic, then $M$ is rigid. More precisely, its universal cover is $\bar N^k\times\Sph^{q}_+(\alpha)$, where $\bar N^k$ is a compact simply connected Einstein manifold with $\Ric_{\bar N}=\lambda g_{\bar N}$.
\end{theorem}

Theorem \ref{thm_criterion} is the counterpart, for quasi-Einstein manifolds, of the theorem of Li, Li, Ou, Qu, Wu and Zhou \cite[Corollary 1.1]{lloqwz}, and we prove it by adapting their argument. As a consequence, Conjecture \ref{conj_main} is reduced to showing that the second fundamental form of $F$ vanishes. To this end, we use that, at the points of $F$, the covariant derivative of the Ricci tensor is determined by $\II$ (Lemma \ref{lem_first_order}). We then derive a Bochner identity for this covariant derivative and integrate it over $F$, which leads to $\II=0$ when $F$ has dimension two or three. Our main result is the following.

\begin{theorem}\label{cor_alldim}
Let $(M^n,g,u,\lambda)$ be a nontrivial compact $m$-quasi-Einstein manifold with boundary, $m>1$ and constant scalar curvature, and suppose that $\operatorname{MAX}(u)$ has dimension $k\in\{2,3\}$. Then $M$ is rigid. More precisely, its universal cover is $\Sph^k(\lambda/(k-1))\times\Sph^{n-k}_+(\alpha)$, and $M$ is the quotient of this product by a finite group of isometries, of order at most two when $k=2$. Consequently, in dimension six, Conjecture \ref{conj_main} holds, except possibly when $k=4$.
\end{theorem}

Here $\Sph^k(\kappa)$ denotes the round sphere of sectional curvature $\kappa$. Theorem \ref{cor_alldim} is the counterpart of \cite[Corollary 1.3]{horfocal}. Some of its cases were already known. Indeed, the cases $k=n-1$ and $(n,k)=(4,2)$ are due to Costa, Ribeiro Jr and Zhou \cite[Theorems 2 and 3]{crz}, and the case $(n,k)=(5,3)$ was considered in \cite{zcao}. To the best of our knowledge, the remaining cases are new, namely $k=2$ in dimension $n\geq5$ and $k=3$ in dimension $n\geq6$. Our argument is the same in all cases, and it does not depend on \cite{zcao}. Since the cases $k=0$ and $k=n-1$ are known, Theorem \ref{cor_alldim} proves Conjecture \ref{conj_main} in dimension five. More precisely, we have the following result.

\begin{theorem}\label{thm_five}
Let $(M^5,g,u,\lambda)$ be a nontrivial compact five-dimensional $m$-quasi-Einstein manifold with boundary, $m>1$ and constant scalar curvature. Then $M$ is rigid. More precisely, if $k$ is the dimension of $\operatorname{MAX}(u)$, then $k\in\{0,2,3,4\}$ and the universal cover of $M$ is
\[
	\Sph^5_+\Big(\frac{\lambda}{m+4}\Big),\quad\ \Sph^2(\lambda)\times\Sph^3_+\Big(\frac{\lambda}{m+2}\Big),\quad\ \Sph^3\Big(\frac{\lambda}{2}\Big)\times\Sph^2_+\Big(\frac{\lambda}{m+1}\Big)\quad\text{or}\quad \widetilde N^4\times\Sph^1_+\Big(\frac{\lambda}{m}\Big)
\]
according as $k=0$, $2$, $3$ or $4$, where $\Sph^k(\kappa)$ is the round sphere of sectional curvature $\kappa$ and $\widetilde N^4$ is a compact simply connected Einstein manifold with $\Ric_{\widetilde N}=\lambda g_{\widetilde N}$. The manifold $M$ is the quotient of its universal cover by a finite group, of order at most two when $k=2$.
\end{theorem}

In Theorem \ref{thm_five}, the cases $k=0$ and $k=4$ are due to He, Petersen and Wylie \cite[Proposition 2.4]{hpw2} and to Costa, Ribeiro Jr and Zhou \cite[Theorem 2]{crz}, respectively, while the case $k=2$, which corresponds to $R=\frac{2(m+5)}{m+2}\lambda$, is new. Note that the bound on the order of the quotient cannot be improved, since it is attained by the product of the real projective plane of curvature $\lambda$ and $\Sph^3_+(\lambda/(m+2))$. Moreover, the example $\mathbb{CP}^2\times\Sph^1_+(\lambda/m)$ shows that the Einstein factor cannot be specified further when $k=4$.

Now suppose that $k=4$. Then $F=\operatorname{MAX}(u)$ is a closed minimal four-manifold. In this case, the integral identity used in the proof of Theorem \ref{cor_alldim} contains a term involving the Weyl tensor $W^F$ of the metric induced on $F$, and this term has no sign. Given an orthonormal frame $\{e_i\}$ of $TF$, we write
\[
	W^F(\II,\II)=\sum_{i,j,\ell,s}W^F(e_i,e_j,e_\ell,e_s)\inner{\II(e_i,e_\ell)}{\II(e_j,e_s)},
\]
and we denote by $D$ the connection induced on $T^*F\otimes T^*F\otimes NF$ by the Levi-Civita connection of $F$ and the normal connection, where $NF$ is the normal bundle of $F$. In our next result, we control this term by means of an integral pinching condition.

\begin{theorem}\label{thm_four}
Let $(M^n,g,u,\lambda)$ be a nontrivial compact $m$-quasi-Einstein manifold with boundary, $m>1$ and constant scalar curvature, and suppose that $F=\operatorname{MAX}(u)$ has dimension $k=4$. Then
\begin{align}\label{eq_four_estimate}
	11\int_FW^F(\II,\II)\dmu_F\geq\int_F\Big(|D\II|^2+\Big(\frac23\lambda+2\alpha\Big)|\II|^2+\frac{1}{72}|\II|^4\Big)\dmu_F,
\end{align}
and, if $n=6$, the coefficient $\frac{1}{72}$ in \eqref{eq_four_estimate} can be replaced by $\frac{11}{8}$. Consequently, if
\begin{align}\label{eq_four_pinching}
	11\int_FW^F(\II,\II)\dmu_F\leq\Big(\frac23\lambda+2\alpha\Big)\int_F|\II|^2\dmu_F,
\end{align}
then $F$ is totally geodesic and $M$ is rigid.
\end{theorem}

In particular, $M$ is rigid whenever $\int_FW^F(\II,\II)\dmu_F\leq0$. Condition \eqref{eq_four_pinching} is satisfied by every rigid manifold with $k=4$, since both sides vanish in this case. Theorem \ref{thm_four} is the counterpart of \cite[Theorem 1.4]{horfocal}. The coefficient of $|\II|^2$ exceeds the one obtained for shrinking solitons by the term $2\alpha$, which tends to zero as $m\to\infty$. When $n=5$, the value $k=4$ already implies rigidity \cite[Theorem 2]{crz}. Therefore the first new case is $n=6$, which is not covered by Theorem \ref{cor_alldim}.

\subsection*{Strategy of proof}

We prove Theorem \ref{thm_criterion} by adapting the argument of \cite{lloqwz}. Let $P=\Ric-\rho g$, where $\rho$ is the eigenvalue of $\Ric$ in the direction of $\nabla u$, which is constant. At the points of $F$, we have $P=m\alpha\,\Pi_T$, where $\Pi_T$ is the orthogonal projection onto $TF$. We work in a tubular neighborhood of $F$. There, we compare $\frac{1}{2m\alpha}P$ with one half of the orthogonal projection onto the orthogonal complement of the fibers of the normal projection onto $F$, and we denote the difference by $E$. Since the first normal derivative of $E$ along $F$ is given by $\II$, the tensor $E$ vanishes to second order when $F$ is totally geodesic. Assume by contradiction that $E$ has a nonzero normal Taylor coefficient, and let $d\geq2$ be the lowest degree of such a coefficient. Then we estimate an averaged curvature term of the normal spheres in two ways. The first estimate comes from the second variation of the total scalar curvature at the round sphere, and the second one from an integral formula of Walczak \cite{walczak}. It turns out that these two estimates are incompatible. The difference with respect to the soliton case is that the shape operator of the level sets of $u$ is $\sqrt{\alpha}\cot(\sqrt{\alpha}\,r)\big(\Id-\frac{1}{m\alpha}P\big)$ instead of $r^{-1}(\Id-\Ric/\lambda)$. We show in Section \ref{sec_criterion} that the additional terms produced by this are of higher order in $r$ (Proposition \ref{prop_higher}). Therefore $E$ vanishes to infinite order along $F$. Since $g$ and $u$ are real analytic in the interior of $M$ \cite[Proposition 2.4]{hpw1}, we conclude that $P^2=m\alpha P$ and $\nabla_{\nabla u}P=0$ on $M$. Rigidity then follows from the theorem of He, Petersen and Wylie on radially flat quasi-Einstein manifolds \cite[Theorem 1.11]{hpw2}.

The proof of Theorem \ref{cor_alldim} takes place at the set $F$. At its points, we have $|\nabla P|^2=3m^2\alpha^2|\II|^2$ (Lemma \ref{lem_first_order}). Since $R$ and $|P|^2$ are constant, the Bochner formula \eqref{eq_bochner_global} for $P$ holds in the interior of $M$. Applying the Laplacian to it, we get a Bochner formula for $\nabla P$ at $F$. The drift term $\frac{m+2}{u}\nabla_{\nabla u}P$ appears in equation \eqref{eq_bochner_P} for $P$. If we compare with \cite{horfocal}, the coefficients depend on the parameter $m$ only through the derivative of this term at $F$ and through $\rho$, the eigenvalue of $\Ric$ in the directions normal to $F$. We then integrate over the closed manifold $F$ and conclude that $\int_FJ\dmu_F=0$, whatever the dimension and the codimension of $F$ (Theorem \ref{thm_general}), where $J$ denotes the integrand of \cite[Theorem 4.6]{horfocal}, computed for $P$ instead of $\Ric$, plus the term $2\alpha|\II|^2$. Since this additional term is nonnegative, the arguments of \cite{horfocal} for focal sets of dimension two and three apply, mutatis mutandis, whence $F$ is totally geodesic when $k=2$ or $k=3$ (Theorems \ref{thm_two} and \ref{thm_three}). Theorems \ref{cor_alldim} and \ref{thm_five} are then immediate consequences of Theorem \ref{thm_criterion}. When $k\geq4$, $J$ contains the pairing $W^F(\II,\II)$, and we prove Theorem \ref{thm_four} as in \cite{horfocal}.

This article is organized as follows. In Section \ref{sec_prelim}, we describe the structure of $u$ near $\operatorname{MAX}(u)$ and prove the Bochner formula for $P$. In Section \ref{sec_bochner}, we derive the Bochner identity at $F$ and its integral form. In Section \ref{sec_lowdim}, we prove that $F$ is totally geodesic when $k\in\{2,3\}$. Section \ref{sec_criterion} is devoted to the proof of Theorem \ref{thm_criterion}, and Section \ref{sec_level} contains the proofs of Theorems \ref{cor_alldim} and \ref{thm_five}. Finally, we prove Theorem \ref{thm_four} in Section \ref{sec_four}.

\subsection*{Acknowledgements}

This study was financed by the Coordena\c{c}\~ao de Aperfei\c{c}oamento de Pessoal de N\'ivel Superior -- Brasil (CAPES) -- Finance Code 001.

\section{Preliminaries}\label{sec_prelim}

Throughout this article, all manifolds are connected. Our convention for the curvature tensor is
\[
	\Rm(X,Y,Z,\cdot)=\inner{\nabla_{Y}\nabla_{X}Z-\nabla_{X}\nabla_{Y}Z+\nabla_{[X,Y]}Z}{\,\cdot\,},
\]
which means that $\Rm(X,Y,X,Y)$ is the sectional curvature of $\mathrm{span}\{X,Y\}$ whenever $X$ and $Y$ are orthonormal, and that $\Ric(Y,Z)=\sum_a\Rm(e_a,Y,e_a,Z)$ for any orthonormal frame $\{e_a\}$. Writing $R_{abcd}=\Rm(e_a,e_b,e_c,e_d)$, the Ricci identity takes the form
\begin{align}\label{eq_ricci_identity}
	\nabla^2_{a,b}\phi_c-\nabla^2_{b,a}\phi_c=\sum_dR_{abcd}\phi_d
\end{align}
for a one-form $\phi$, and the same holds for each entry of a tensor. The contracted second Bianchi identity is $\sum_a\nabla_aR_{abcd}=\nabla_c\Ric_{bd}-\nabla_d\Ric_{bc}$. A symmetric bilinear form and the corresponding self-adjoint endomorphism are denoted by the same symbol. Unless otherwise stated, the norm of a tensor is its full tensor norm. Given symmetric bilinear forms $\mathcal{A}$ and $\mathcal{B}$, we write
\[
	\Rm(\mathcal{A})_{bd}=\sum_{a,c}R_{abcd}\mathcal{A}_{ac},\qquad \Rm(\mathcal{A},\mathcal{B})=\sum_{a,b,c,d}R_{abcd}\mathcal{A}_{ac}\mathcal{B}_{bd}.
\]

From now on, we assume that $(M^n,g,u,\lambda)$ is a nontrivial compact $m$-quasi-Einstein manifold with nonempty boundary, $m>1$ and constant scalar curvature $R$. We write $U=\max_Mu$, and we have $\lambda>0$ by \cite[Theorem 4.1]{hpw1}. Following \cite{hpw2}, we set
\begin{align}\label{eq_rho_P}
	\rho=\frac{(n-1)\lambda-R}{m-1},\qquad P=\Ric-\rho g,\qquad \alpha=\frac{\lambda-\rho}{m},
\end{align}
so that equation \eqref{fundeq} becomes $\nabla^2u=\frac{u}{m}(P-m\alpha g)$. Since $R$ is constant, it follows from \cite[Proposition 3.3 and Corollary 3.5]{hpw2} that
\begin{align}\label{eq_P_identities}
	P(\nabla u,\cdot)=0,\qquad \divg P=0,\qquad |P|^2=m\alpha\,\trg P,\qquad \alpha>0.
\end{align}
In particular, $\nabla u$ is an eigenvector of $\Ric$ with eigenvalue $\rho$, and $\trg P=R-n\rho$ is constant. Using the first identity in \eqref{eq_P_identities}, we get $\nabla|\nabla u|^2=2\nabla^2u(\nabla u,\cdot)=-2\alpha u\nabla u$. Thus $|\nabla u|^2+\alpha u^2$ is constant on $M$, and evaluating it at a maximum point of $u$ we conclude that
\begin{align}\label{eq_transnormal}
	|\nabla u|^2=\alpha\big(U^2-u^2\big),
\end{align}
and, consequently, the critical points of $u$ are exactly its maximum points; see also \cite[Proposition 3.11]{hpw2} and \cite[Remark 4]{crz}. It is convenient to introduce the function
\begin{align}\label{eq_ab}
	\psi(r)=\frac{\sin(\sqrt{\alpha}\,r)}{\sqrt{\alpha}},\qquad\text{so that}\qquad \psi''=-\alpha\psi,\qquad (\psi')^2=1-\alpha\psi^2,\qquad \Big(\frac{\psi'}{\psi}\Big)'=-\frac{1}{\psi^2}.
\end{align}

\begin{proposition}\label{prop_structure}
The set $F=\operatorname{MAX}(u)$ has the following properties.
\begin{enumerate}[label=\textup{(\roman*)},leftmargin=*]
\item $F$ is a closed connected embedded submanifold contained in the interior of $M$, whose dimension we denote by $k$. At the points of $F$, $\nabla^2u=-\alpha U\,\Pi_N$, so that $P=m\alpha\,\Pi_T$ and $\Ric=\lambda\Pi_T+\rho\Pi_N$, where $\Pi_T$ and $\Pi_N$ are the orthogonal projections onto $TF$ and onto the normal bundle $NF$ of $F$.
\item With $q=n-k$, we have
\begin{align}\label{eq_constants}
	\rho=(q-1)\alpha,\qquad \lambda=(m+q-1)\alpha,\qquad \trg P=km\alpha,\qquad |P|^2=km^2\alpha^2,
\end{align}
and $R=k\lambda+q\rho$ is given by \eqref{eq_quantization}.
\item The distance $r$ to $F$ satisfies
\begin{align}\label{eq_u_cos}
	u=U\cos(\sqrt{\alpha}\,r)\qquad\text{on } M,
\end{align}
and there is $\varepsilon>0$ such that the set $F_\varepsilon$ of the points of $M$ at distance less than $\varepsilon$ from $F$ is contained in the interior of $M$ and the normal exponential map of $F$ is a diffeomorphism from $\{\xi\in NF:\ |\xi|<\varepsilon\}$ onto $F_\varepsilon$.
\item On $F_\varepsilon\setminus F$, let $\nu=\nabla r$, let $\Sigma_r$ be the level set of $r$ through a given point, and let $SX=\nabla_X\nu$ be its shape operator. Then $P(\nu,\cdot)=0$ and, on $T\Sigma_r$,
\begin{align}\label{eq_shape}
	S=\frac{\psi'}{\psi}\Big(\Id-\frac{1}{m\alpha}P\Big),\qquad \trg S=(q-1)\frac{\psi'}{\psi},\qquad |S|^2=(q-1)\Big(\frac{\psi'}{\psi}\Big)^2.
\end{align}
Moreover, the normal exponential map sends the bundle of normal vectors of length $r$ onto $\Sigma_r$ and, under the identification with the unit normal bundle of $F$, the volume measure of $\Sigma_r$ is $\psi(r)^{q-1}$ times the product of the measure of $F$ and of the standard measure of the unit spheres of $NF$.
\end{enumerate}
\end{proposition}

\begin{proof}
(i) We set $r=\alpha^{-1/2}\arccos(u/U)$. This function is continuous on $M$, and its zero set is $F$. Moreover, it follows from \eqref{eq_transnormal} that $r$ is smooth on $M\setminus F$, where $|\nabla r|=1$. Our first claim is that $r$ is the distance to $F$. Indeed, consider a curve joining $x$ to a point of $F$. Since $|\nabla r|=1$, the length of this curve is at least $r(x)$, and consequently $r$ does not exceed the distance to $F$. To prove the reverse inequality, let us arbitrarily fix $x\in M\setminus F$ and consider the integral curve of $-\nabla r$ starting at $x$, which is a unit speed geodesic because $\nabla_{\nabla r}\nabla r=0$. Since $u$ increases along this curve, the curve does not leave $M$, and $r$ decreases at unit rate along it. By the compactness of $M$, after length $r(x)$ the curve converges to a point of $F$, which we will denote by $\pi(x)$. Therefore $r$ is the distance to $F$, as claimed. Setting $\pi(x)=x$ for $x\in F$, we have a map $\pi\colon M\to F$, and we will now check that it is continuous. If $x\in M\setminus F$, then $\pi(x)$ is the point at distance $r(x)$ from $x$ on the geodesic with initial velocity $-\nabla r(x)$. Thus $\pi$ is continuous on $M\setminus F$, by the continuous dependence of geodesics on their initial data. On the other hand, given $y\in F$, the distance from $\pi(x)$ to $y$ is at most $r(x)$ plus the distance from $x$ to $y$, which is at most twice the distance from $x$ to $y$. A fortiori, $\pi$ is continuous at the points of $F$ as well. Since $M$ is connected and $\pi(M)=F$, it then follows that $F$ is connected. By the argument of \cite[Lemma 7]{wang} for transnormal functions, $F$ is a smooth submanifold; see the proof of \cite[Theorem 1]{crz}. Furthermore, by \cite[Proposition 3.13]{hpw2}, the vectors tangent to $F$ are eigenvectors of $P$ with eigenvalue $m\alpha$, while the vectors normal to $F$ are eigenvectors of $P$ with eigenvalue $0$. Since $u=0$ on $\partial M$, $F$ is contained in the interior of $M$. It is compact, because it is closed in $M$. Finally, the formula for $\nabla^2u$ follows from $\nabla^2u=\frac{u}{m}(P-m\alpha g)$.

(ii) Since $\trg P$ is constant, we can compute it at a point of $F$, which yields $\trg P=km\alpha$. Then \eqref{eq_P_identities} implies that $|P|^2=km^2\alpha^2$. Taking the trace of $\Ric=\lambda\Pi_T+\rho\Pi_N$, we get $R=k\lambda+q\rho$. Inserting this into the definition of $\rho$, a direct computation shows that $(m+q-1)\rho=(q-1)\lambda$, whence \eqref{eq_constants} and \eqref{eq_quantization} follow from \eqref{eq_rho_P}.

(iii) The identity \eqref{eq_u_cos} was established in the proof of (i). Since $F$ is compact, there is $\varepsilon\in(0,\pi/(2\sqrt\alpha))$ such that the normal exponential map of $F$ is a diffeomorphism from $\{\xi\in NF:\ |\xi|<\varepsilon\}$ onto $F_\varepsilon$. For such $\varepsilon$, we have $u>0$ on $F_\varepsilon$, that is, $F_\varepsilon$ is contained in the interior of $M$. By the proof of (i), the integral curves of $\nabla r$ on $F_\varepsilon\setminus F$ are unit speed geodesics which realize the distance to $F$. A fortiori, they are the normal geodesics of $F$.

(iv) Since $\nu$ is a multiple of $\nabla u$, we have $P(\nu,\cdot)=0$ by \eqref{eq_P_identities}. It follows from \eqref{eq_u_cos} that $u=U\psi'(r)$, and consequently $\nabla^2u(X,Y)=-\alpha U\psi\inner{SX}{Y}$ for $X,Y$ tangent to $\Sigma_r$. Comparing this with $\nabla^2u=\frac{u}{m}(P-m\alpha g)$, we get the first identity in \eqref{eq_shape}. The other two follow from \eqref{eq_constants} and $P(\nu,\cdot)=0$. Finally, note that the density of the volume measure of $\Sigma_r$ with respect to the product measure has logarithmic derivative $\trg S=(q-1)\psi'/\psi$ with respect to $r$, and that it is asymptotic to $r^{q-1}$ as $r\to0$. Therefore it is equal to $\psi^{q-1}$, as desired.
\end{proof}

We have $k\neq1$ by \cite[Proposition 5]{crz}. If $k=0$, then $P=0$ by \eqref{eq_constants}, and $M$ is Einstein. In what follows, $F=\operatorname{MAX}(u)$, $g_F$ is the metric induced on $F$, $|F|$ is its volume, $k$ is its dimension and $q=n-k$, while $\varepsilon$ and $F_\varepsilon$ are as in Proposition \ref{prop_structure}. On $F_\varepsilon\setminus F$, the integral curves of $\nu$ are the normal geodesics of $F$. In order to identify tensors on different level sets, we use parallel orthonormal frames along these geodesics, tangent to the level sets. Given a tensor on the level sets, a prime denotes its covariant derivative in the direction of $\nu$, which coincides with the derivative of its matrix in such a frame. We denote by $\Ric_\Sigma$ the Ricci operator of the metric induced on $\Sigma_r$. The unit geodesic field $\nu$ satisfies the Riccati equation
\begin{align}\label{eq_riccati}
	\Rm(\nu,X,\nu,Y)=-\inner{(S'+S^2)X}{Y},\qquad X,Y\perp\nu,
\end{align}
and the trace of the Gauss equation is $\Ric_\Sigma=\Ric-\Rm(\nu,\cdot,\nu,\cdot)+(\trg S)S-S^2$ on $T\Sigma_r$. Combining these two equations, we have $\Ric_\Sigma=\Ric+S'+(\trg S)S$. We now insert \eqref{eq_shape} and $\Ric=P+\rho g$ into this identity. Using \eqref{eq_ab} and \eqref{eq_constants}, we conclude that, on $T\Sigma_r$,
\begin{align}\label{eq_ric_sigma}
	\Ric_\Sigma=\frac{q-2}{\psi^2}\Id+\frac{1}{m\alpha}\Big(\lambda-\frac{q-2}{\psi^2}\Big)P-\frac{\psi'}{m\alpha\psi}\nabla_\nu P.
\end{align}

In the next lemma, we prove the Bochner formula for $P$. It is worth mentioning that, when $q=2$, the second identity below is \cite[Corollary 3.2]{zcao}; see also \cite[Lemma 4]{crz}.

\begin{lemma}\label{lem_bochner_P}
In the interior of $M$, let $f=-m\log u$. Then
\begin{align}
	(\nabla_XP)(Y,Z)-(\nabla_YP)(X,Z)&=\Rm(X,Y,\nabla f,Z)\nonumber\\
	&\quad+\frac1m\big(\mathrm{d}f(Y)(m\alpha g-P)(X,Z)-\mathrm{d}f(X)(m\alpha g-P)(Y,Z)\big),\label{eq_codazzi_P}\\
	\Delta P+\frac{m+2}{u}\nabla_{\nabla u}P&=2\lambda P-2\Rm(P),\label{eq_bochner_P}\\
	0&=|\nabla P|^2+2\lambda|P|^2-2\Rm(P,P).\label{eq_bochner_global}
\end{align}
\end{lemma}

The proof of Lemma \ref{lem_bochner_P} is a computation in an orthonormal frame, and we give it in Appendix \ref{app_proofs}.

Most of the computations below are carried out at a point of $F$, and we fix the notation for them here. Let us arbitrarily fix $p\in F$ and orthonormal bases $\{e_i\}$ of $T_pF$ and $\{e_\beta\}$ of $N_pF$. The indices $i,j,\ell,s$ are tangent, the indices $\beta,\gamma,\delta,\eta$ are normal, and the indices $a,b,c,d,x,y,z$ run over all directions. Repeated indices are summed. The second fundamental form of $F$ is $\II(X,Y)=(\nabla_XY)^\perp$, with components $h^\beta_{ij}=\inner{\II(e_i,e_j)}{e_\beta}$. Given a normal vector $\xi$, we denote by $S_\xi$ the shape operator in the direction of $\xi$, that is, $\inner{S_\xi X}{Y}=\inner{\II(X,Y)}{\xi}$, and we write $S_\beta$ when $\xi=e_\beta$. We denote by $D$ the connection induced on $T^*F\otimes T^*F\otimes NF$ by the Levi-Civita connection of $F$ and the normal connection $\nabla^\perp$, and by $R^\perp(X,Y)=\nabla^\perp_X\nabla^\perp_Y-\nabla^\perp_Y\nabla^\perp_X-\nabla^\perp_{[X,Y]}$ the curvature of the normal connection. Products of shape operators are understood as products of matrices in the basis $\{e_i\}$, that is, $(S_\beta S_\gamma)_{ij}=\sum_\ell h^\beta_{i\ell}h^\gamma_{\ell j}$. We shall also use the following quantities at $p$:
\begin{center}
\renewcommand{\arraystretch}{1.6}
\begin{tabular}{@{}l@{\qquad}l@{}}
	shape operators & $t=|\II|^2$,\quad $G_{\beta\gamma}=\trg(S_\beta S_\gamma)$,\quad $c=|G|^2$,\quad $T=\sum_\beta S_\beta^2$,\\
	& $\mathcal{P}_{\beta\gamma}=S_\beta S_\gamma+S_\gamma S_\beta$,\quad $K_{\beta\gamma}=S_\beta S_\gamma-S_\gamma S_\beta$,\\
	& $\sigma_{ij\ell s}=\inner{\II_{ij}}{\II_{\ell s}}+\inner{\II_{i\ell}}{\II_{js}}+\inner{\II_{is}}{\II_{j\ell}}$,\\
	divergence of $\II$ & $d^\beta_j=\sum_i(D_{e_i}h^\beta)_{ij}$,\\
	curvature of $M$ & $L_{\beta\gamma ij}=R_{\beta i\gamma j}$,\quad $V_{\beta\gamma ij}=R_{\beta\gamma ij}$,\quad $Y_{\beta i\gamma\delta}=R_{\beta i\gamma\delta}$,\\
	& $\mathcal{R}_{ij}=\sum_\ell R_{\ell i\ell j}$,\quad $\tilde{\mathcal{R}}=\Rm(\Pi_T)$, so that $\tilde{\mathcal{R}}_{bd}=\sum_iR_{ibid}$.
\end{tabular}
\end{center}
The tensor $\sigma$ is totally symmetric on $T_pF$, and $L^s$ denotes the symmetrization of $L$ in $(i,j)$. When we write a pairing such as $\inner{K}{V}=\sum(K_{\beta\gamma})_{ij}V_{\beta\gamma ij}$ or $\inner{\mathcal{P}}{L^s}=\sum(\mathcal{P}_{\beta\gamma})_{ij}L^s_{\beta\gamma ij}$, all indices are contracted.

\needspace{4\baselineskip}
\begin{lemma}\label{lem_first_order}
At $p\in F$:
\begin{enumerate}[label=\textup{(\roman*)},leftmargin=*]
\item $\nabla P$ is totally symmetric; its only nonzero components are those with exactly one normal entry, and $(\nabla P)(e_i,e_j,e_\beta)=m\alpha h^\beta_{ij}$. In particular, $|\nabla P|^2=3m^2\alpha^2t$, and $F$ is minimal.
\item (Codazzi equation) $R_{i\ell j\beta}=(D_\ell h^\beta)_{ij}-(D_ih^\beta)_{\ell j}$, and $\sum_iR_{i\beta ij}=-d^\beta_j$.
\item (Ricci equation) $\inner{R^\perp(e_i,e_j)e_\gamma}{e_\beta}=V_{\beta\gamma ij}+(K_{\beta\gamma})_{ij}$.
\item For all $a,b,c,z$,
\[
	\nabla^2_{z,a}P_{bc}-\nabla^2_{z,b}P_{ac}=m\alpha\sum_dR_{abdc}(\Pi_N)_{zd}+m\alpha^2\big((\Pi_N)_{zb}(\Pi_N)_{ac}-(\Pi_N)_{za}(\Pi_N)_{bc}\big).
\]
\end{enumerate}
\end{lemma}

\begin{proof}
(i) Since $\nabla f(p)=0$, it follows from \eqref{eq_codazzi_P} that, at $p$, the tensor $\nabla P$ is symmetric in its first two entries. A fortiori, it is totally symmetric. Let $X,Y,Z$ be tangent and $\xi,\zeta$ be normal to $F$. If we differentiate $P=m\alpha\,\Pi_T$ along $F$, then
\[
	(\nabla_XP)(Y,Z)=0,\qquad (\nabla_XP)(\xi,\zeta)=0,\qquad (\nabla_XP)(Y,\xi)=m\alpha\inner{\II(X,Y)}{\xi},
\]
because the tangent part of $\nabla_X\xi$ is $-S_\xi X$. The components with a tangent entry are then determined by total symmetry. As for the remaining ones, let us arbitrarily fix a unit normal vector $\xi$ at $p$. By Proposition \ref{prop_structure}, the velocity of the geodesic $r\mapsto\exp_p(r\xi)$ is $\nu$ for $r>0$, whence $P$ vanishes on this velocity. Taking the derivative at $r=0$, we conclude that $(\nabla_\xi P)(\xi,\xi)=0$. Now the restriction of $\nabla P$ to $N_pF$ is a totally symmetric cubic form, and it then follows by polarization that the components with three normal entries vanish. Since the normal entry can occupy three positions, we have $|\nabla P|^2=3m^2\alpha^2t$. Finally, $\trg P$ is constant, and therefore $0=\sum_i(\nabla_\beta P)_{ii}=m\alpha\sum_ih^\beta_{ii}$, that is, $F$ is minimal.

(ii) The first identity is the Codazzi equation, which is the normal part of the definition of $\Rm$. For the second one, we sum over $i=j$ and use the minimality of $F$, which leads to $\sum_iR_{i\ell i\beta}=-d^\beta_\ell$. The assertion then follows from $R_{i\beta i\ell}=R_{i\ell i\beta}$.

(iii) Since $\nabla_X\xi=\nabla^\perp_X\xi-S_\xi X$ for a normal field $\xi$, the normal part of the definition of $\Rm$ is
\[
	\Rm(X,Y,\xi,\zeta)=-\inner{R^\perp(X,Y)\xi}{\zeta}+\inner{[S_\xi,S_\zeta]X}{Y}.
\]
Taking $X=e_i$, $Y=e_j$, $\xi=e_\beta$ and $\zeta=e_\gamma$, and recalling the matrix convention fixed above, we have $\inner{[S_\beta,S_\gamma]e_i}{e_j}=-(K_{\beta\gamma})_{ij}$, as claimed.

(iv) Immediate consequence of differentiating \eqref{eq_codazzi_P} in the direction of $e_z$, recalling that $\nabla f=0$, $\nabla^2f=m\alpha\,\Pi_N$ and $m\alpha g-P=m\alpha\,\Pi_N$ at $p$, by Proposition \ref{prop_structure}(i).
\end{proof}

\section{A Bochner identity at the maximum set}\label{sec_bochner}

The main goal of this section is to prove the identity $\int_FJ\dmu_F=0$ of Theorem \ref{thm_general}, on which Section \ref{sec_lowdim} is based. For shrinking solitons, this identity is \cite[Theorem 4.6]{horfocal}, and it is proved in \cite[Sections 3 and 4]{horfocal} by a computation at the points of the focal set. In that computation, the roles of $P$ and $m\alpha$ are played by $\Ric$ and $\lambda$, and the soliton equation is used only through the counterparts of the following facts: the identity $\nabla\Ric=\nabla P$, the values of $P$, $\nabla P$ and $\Ric$ at the points of $F$ (Proposition \ref{prop_structure}(i) and Lemma \ref{lem_first_order}(i)), the identity of Lemma \ref{lem_first_order}(iv), the constancy of $|P|^2$ and the Bochner formulas \eqref{eq_bochner_P} and \eqref{eq_bochner_global}. These facts differ from their counterparts in three places. The identity of Lemma \ref{lem_first_order}(iv) has an additional term, equation \eqref{eq_bochner_P} has a different drift term, and the eigenvalue $\rho$ of $\Ric$ in the directions normal to $F$ is zero for solitons. Hence several statements of \cite{horfocal} hold here with the same proof, and for these we will only indicate what the proof uses. We will give the details where the parameter $m$ enters. Our notation is that of \cite{horfocal}, with the following changes. The focal set, its normal bundle and the tensors denoted in \cite{horfocal} by $P_{\alpha\beta}$, $\rho$, $\tilde\rho$ and $F$ are denoted here by $F$, $NF$, $\mathcal{P}_{\beta\gamma}$, $\mathcal{R}$, $\tilde{\mathcal{R}}$ and $\mathcal{F}$. The normal indices are $\beta,\gamma,\delta,\eta$, and the indices $m$ and $e$ of \cite{horfocal} are replaced by $s$ and $z$.

At $p$, we write $\Lambda_{zabc}=(\nabla^2_{z,a}P)(e_b,e_c)$ and $\mathcal{F}_{\beta\gamma\delta\eta}=\Lambda_{\beta\gamma\delta\eta}$, and we denote the components of $D\II$ by $h^\beta_{ij,\ell}=(D_\ell h^\beta)_{ij}$.

\begin{lemma}\label{lem_hessian}
At $p\in F$,
\begin{align}\label{eq_hessian_norm}
	\frac{1}{m^2\alpha^2}|\nabla^2P|^2=|\sigma|^2+4|D\II|^2+6|\mathcal{P}|^2+8\inner{\mathcal{P}}{L^s}+6|L^s|^2+\frac12|V|^2+2|Y|^2+\frac{1}{m^2\alpha^2}|\mathcal{F}|^2.
\end{align}
\end{lemma}

\begin{proof}
In \cite[Lemma 3.1]{horfocal}, every component of the Hessian of the Ricci tensor with at least one tangent entry is computed, and the norm is obtained by adding the squares. That proof uses the components of the first derivative along the focal set, which are given here by Lemma \ref{lem_first_order}(i), the Ricci identity \eqref{eq_ricci_identity}, and the identity of Lemma \ref{lem_first_order}(iv) in two cases, namely when $z$ is tangent, and when $z$ is normal and one of the indices $a,b,c$ is tangent. Note that, in both cases, the last term of Lemma \ref{lem_first_order}(iv) vanishes. Mutatis mutandis, the same proof then shows that the components of $\nabla^2P$ with at least one tangent entry are those of \cite[Lemma 3.1]{horfocal}, with $\lambda$ replaced by $m\alpha$, whence \eqref{eq_hessian_norm} follows. The components $\mathcal{F}$ are not computed, and they enter \eqref{eq_hessian_norm} only through $|\mathcal{F}|^2$.
\end{proof}

\begin{lemma}\label{lem_hess_two}
At $p\in F$:
\begin{enumerate}[label=\textup{(\roman*)},leftmargin=*]
\item $\sum_iR_{i\beta i\gamma}=-\frac32G_{\beta\gamma}$ and $\sum_\delta R_{\delta\beta\delta\gamma}=\rho\,\delta_{\beta\gamma}+\frac32G_{\beta\gamma}$;
\item $\mathcal{R}=\Ric_F+T$, $\trg\mathcal{R}=k\lambda+\frac32t$, and $|\tilde{\mathcal{R}}|^2=|\mathcal{R}|^2+2|d|^2+\frac94c$.
\end{enumerate}
In particular, the scalar curvature of $F$ is $R_F=k\lambda+\frac12t$. If $k=2$, then $\mathcal{R}=(\lambda+\frac34t)\Id$ and the Gaussian curvature of $F$ is $\lambda+\frac14t$.
\end{lemma}

\begin{proof}
The first identity in (i) and the assertions in (ii) are proved as in \cite[Lemma 3.2]{horfocal}. Indeed, that proof uses the second derivative of the constant $|P|^2$ in the normal directions at $p$, the components of $\nabla^2P$ computed in \cite[Lemma 3.1]{horfocal}, the Gauss equation, the minimality of $F$, Lemma \ref{lem_first_order}(ii) and the identity $\Ric(e_i,e_i)=\lambda$. And for the second identity in (i), it suffices to note that $\sum_iR_{i\beta i\gamma}+\sum_\delta R_{\delta\beta\delta\gamma}=\Ric(e_\beta,e_\gamma)=\rho\,\delta_{\beta\gamma}$, by Proposition \ref{prop_structure}(i). Here the eigenvalue $\rho$ appears, which is zero in \cite{horfocal}.
\end{proof}

We will now apply the Laplacian to the identity \eqref{eq_bochner_global}, which holds in the interior of $M$, and evaluate the result at the arbitrarily fixed point $p\in F$. Since $|P|^2$ is constant, a direct computation shows that
\begin{align}\label{eq_bochner_laplacian}
	|\nabla^2P|^2+\inner{\nabla P}{\Delta\nabla P}=(\Delta\Rm)(P,P)+2\Rm(\Delta P,P)+2\Theta+4m^2\alpha^2\Omega,
\end{align}
where $\Delta\Rm$ is the rough Laplacian of the curvature tensor and
\[
	\Theta=\sum_z\Rm(\nabla_zP,\nabla_zP),\qquad \Omega=\frac{1}{m\alpha}\sum_{z,i,b,d}(\nabla_z\Rm)_{ibid}(\nabla_zP)_{bd},
\]
with $i$ tangent and $z,b,d$ arbitrary. The last two terms of \eqref{eq_bochner_laplacian} are the products of first derivatives $\sum_z\big(2\Rm(\nabla_zP,\nabla_zP)+4(\nabla_z\Rm)(P,\nabla_zP)\big)$, written with $P=m\alpha\,\Pi_T$ at $p$. The first term on the left-hand side is given by Lemma \ref{lem_hessian}, and we will now compute the other three.

We begin with $2\Rm(\Delta P,P)$. Since $\nabla u=0$ at $p$, equation \eqref{eq_bochner_P} reduces to $\Delta P=2\lambda P-2\Rm(P)$ at $p$. Moreover, $P=m\alpha\,\Pi_T$ implies that $\Rm(P)=m\alpha\tilde{\mathcal{R}}$ and that $\Rm(\mathcal{A},P)=m\alpha\sum_{a,c}\tilde{\mathcal{R}}_{ac}\mathcal{A}_{ac}$ for every symmetric $\mathcal{A}$. Therefore, at $p$,
\begin{align}
	\Delta P&=2\lambda m\alpha\,\Pi_T-2m\alpha\tilde{\mathcal{R}},\label{eq_lap_ric}\\
	2\Rm(\Delta P,P)&=4\lambda m^2\alpha^2\trg\mathcal{R}-4m^2\alpha^2|\tilde{\mathcal{R}}|^2.\label{eq_term_ricci}
\end{align}

\begin{lemma}\label{lem_term_rm}
At $p\in F$,
\[
	(\Delta\Rm)(P,P)=m^2\alpha^2\Big(2\lambda\trg\mathcal{R}-2|\tilde{\mathcal{R}}|^2+2|L^s|^2-\frac32|V|^2\Big).
\]
\end{lemma}

\begin{proof}
Mutatis mutandis, this is \cite[Lemma 4.1]{horfocal}. Indeed, the proof given there uses the second Bianchi identity, the Ricci identity \eqref{eq_ricci_identity}, the components $\Lambda_{\ell ijs}=-m\alpha\sigma_{\ell ijs}$ of $\nabla^2\Ric=\nabla^2P$ with four tangent entries, and the identities $\Ric(e_i,e_j)=\lambda\delta_{ij}$ and $\Ric(e_i,e_\beta)=0$, which hold here by Proposition \ref{prop_structure}(i).
\end{proof}

\begin{lemma}\label{lem_term_grad}
At $p\in F$,
\[
	\inner{\nabla P}{\Delta\nabla P}=(9\lambda+2\alpha)m^2\alpha^2t-2m^2\alpha^2\Omega-6\Theta.
\]
\end{lemma}

For solitons, the corresponding identity is \cite[Lemma 4.2]{horfocal}, in which the coefficient of $t$ is $9\lambda^3$. Its proof commutes the Laplacian with the covariant derivative, and it is easy to see that only two terms of that computation change here. The first is $\Ric_{zx}\nabla_xP_{ab}$, in which $\Ric_{zx}=\lambda(\Pi_T)_{zx}+\rho(\Pi_N)_{zx}$. The second is the derivative of the drift term of \eqref{eq_bochner_P}, which equals $(m+2)\alpha(\Pi_N)_{zx}\nabla_xP_{ab}$ at $p$. Since $\rho+(m+2)\alpha=\lambda+2\alpha$ by \eqref{eq_constants}, their sum is $\lambda\nabla_zP_{ab}+2\alpha(\Pi_N)_{zx}\nabla_xP_{ab}$, while for solitons it is $\lambda\nabla_z\Ric_{ab}$. After contraction with $\nabla_zP_{ab}$, the additional term gives $2\alpha\sum_\beta|\nabla_\beta P|^2=2\alpha m^2\alpha^2t$, as claimed. For the convenience of the reader, the complete computation is given in Appendix \ref{app_proofs}.

\begin{proposition}\label{prop_bochner}
At every point of $F$,
\begin{align}\label{eq_bochner_point}
	\frac{1}{m^2\alpha^2}|\nabla^2P|^2+6|\mathcal{R}-\lambda\Id|^2+(18\lambda+2\alpha)t+12|d|^2+\frac{27}{2}c+\frac32|V|^2-2|L^s|^2=\frac{8}{m^2\alpha^2}\Theta+6\Omega.
\end{align}
\end{proposition}

\begin{proof}
Inserting \eqref{eq_term_ricci} and the identities of Lemmas \ref{lem_term_rm} and \ref{lem_term_grad} into \eqref{eq_bochner_laplacian}, we have
\begin{align*}
	|\nabla^2P|^2&+(9\lambda+2\alpha)m^2\alpha^2t-2m^2\alpha^2\Omega-6\Theta\\
	&=m^2\alpha^2\Big(6\lambda\trg\mathcal{R}-6|\tilde{\mathcal{R}}|^2+2|L^s|^2-\frac32|V|^2\Big)+2\Theta+4m^2\alpha^2\Omega.
\end{align*}
We divide this identity by $m^2\alpha^2$. Since Lemma \ref{lem_hess_two}(ii) implies that
\[
	6\lambda\trg\mathcal{R}-6|\tilde{\mathcal{R}}|^2-9\lambda t=-6|\mathcal{R}-\lambda\Id|^2-18\lambda t-12|d|^2-\frac{27}{2}c,
\]
it then follows that \eqref{eq_bochner_point} holds, as desired.
\end{proof}

It remains to write $\Theta$ and $\Omega$ in terms of the geometry of $F$. Let us set
\[
	\omega=\sum_\beta\sum_{i,j,\ell,s}R_{ij\ell s}h^\beta_{i\ell}h^\beta_{js},\qquad A^\beta_{\ell ij}=h^\beta_{ij,\ell}-h^\beta_{\ell j,i}.
\]
The tensor $A^\beta$ is skew in $(\ell,i)$, and, with $C_D=\sum_{\beta,i,j,\ell}(D_ih^\beta)_{j\ell}(D_jh^\beta)_{i\ell}$, we have $|A|^2=2|D\II|^2-2C_D$.

\begin{lemma}\label{lem_X}
At $p\in F$, $(m\alpha)^{-2}\Theta=\omega+\frac32\inner{K}{V}-\inner{\mathcal{P}}{L^s}$. Moreover,
\[
	\int_F\Omega\dmu_F=\int_F\Big(3|d|^2-\frac12|A|^2+3\inner{T}{\mathcal{R}}+\omega+2\inner{\mathcal{P}}{L^s}-\frac32\inner{K}{V}+\frac92c\Big)\dmu_F.
\]
\end{lemma}

\begin{proof}
The formula for $\Theta$ is \cite[Lemma 4.4]{horfocal}, with $\lambda$ replaced by $m\alpha$, and its proof uses only the components of $\nabla P$ at $p$ given by Lemma \ref{lem_first_order}(i). By \cite[Lemma 4.5]{horfocal}, $\Omega$ is the sum of $3\tau_3-\tau_1$ and of the last five terms of the integrand above, where $\tau_3$ and $\tau_1$ are the functions on $F$ defined in \cite[Section 4]{horfocal}, which satisfy $\int_F\tau_3\dmu_F=\int_F|d|^2\dmu_F$ and $\int_F\tau_1\dmu_F=\int_F(|D\II|^2-C_D)\dmu_F$. The proof of that lemma, including the computation of the connection terms in \cite[Appendix A]{horfocal}, uses Lemma \ref{lem_first_order}(i) and (ii), the second Bianchi identity, the minimality of $F$ and the first identity of Lemma \ref{lem_hess_two}(i). Recalling that $|D\II|^2-C_D=\frac12|A|^2$, we see that the formula for the integral of $\Omega$ follows.
\end{proof}

For later use, we record the inequality
\begin{align}\label{eq_d_bound}
	|d|^2\leq\frac{k-1}{2}|A|^2,
\end{align}
with equality when $k=2$, and the identities
\begin{align}\label{eq_algebraic_norms}
\begin{split}
	\sum_{\beta,\gamma}\trg(S_\beta S_\gamma S_\beta S_\gamma)&=|T|^2-C,\qquad |\mathcal{P}|^2=4|T|^2-2C,\\
	|K|^2&=2C,\qquad |\sigma|^2=3c+6|T|^2-6C.
\end{split}
\end{align}
where $C=\frac12\sum_{\beta,\gamma}|K_{\beta\gamma}|^2$; see \cite[Section 4]{horfocal}.

\begin{theorem}\label{thm_general}
For every $k\geq2$ and $q\geq1$, $\int_FJ\dmu_F=0$, where
\begin{align*}
	J={}&4|D\II|^2+3|A|^2-6|d|^2+4\Big|L^s+\frac12\mathcal{P}\Big|^2+2\Big|V-\frac34K\Big|^2+2|Y|^2+\frac{1}{m^2\alpha^2}|\mathcal{F}|^2\\
	&+6|\mathcal{R}-\lambda\Id|^2+(18\lambda+2\alpha)t+|\sigma|^2+5|\mathcal{P}|^2-\frac94C-\frac{27}{2}c-14\omega-18\inner{T}{\mathcal{R}}.
\end{align*}
\end{theorem}

\begin{proof}
We insert \eqref{eq_hessian_norm} and the formula for $\Theta$ of Lemma \ref{lem_X} into \eqref{eq_bochner_point}, we integrate over $F$, and we use Lemma \ref{lem_X} for the integral of $\Omega$. The terms with derivatives of $\II$ add up to $4|D\II|^2+3|A|^2-6|d|^2$. In the remaining terms, we complete the squares by means of $6|\mathcal{P}|^2+4\inner{\mathcal{P}}{L^s}+4|L^s|^2=4|L^s+\frac12\mathcal{P}|^2+5|\mathcal{P}|^2$ and
\[
	2|V|^2-3\inner{K}{V}=2\Big|V-\frac34K\Big|^2-\frac98|K|^2=2\Big|V-\frac34K\Big|^2-\frac94C,
\]
where the last equality follows from \eqref{eq_algebraic_norms}. This gives $J$, as claimed. Mutatis mutandis, the computation is the one in the proof of \cite[Theorem 4.6]{horfocal}, the only difference being the term $2\alpha t$, which comes from \eqref{eq_bochner_point}.
\end{proof}

In Sections \ref{sec_lowdim} and \ref{sec_four}, we use $J$ in terms of the intrinsic geometry of $F$. By Lemma \ref{lem_hess_two}(ii), the trace free part of $\Ric_F$ is $\mathring{\Ric}_F=\Ric_F-(\lambda+\frac{t}{2k})g_F$. We denote by $\mathring{T}=T-\frac{t}{k}\Id$ the trace free part of $T$, and we write $R^F(S,S)=\sum R^F_{ij\ell s}S_{i\ell}S_{js}$ and $W^F(\II,\II)=\sum_\beta W^F(S_\beta,S_\beta)$, where $W^F(S,S)$ is defined in the same way. The computation which follows \cite[Theorem 4.6]{horfocal} uses Lemma \ref{lem_hess_two}(ii), the Gauss equation, the decomposition of the curvature tensor of $F$ into its Weyl and Schouten parts and \eqref{eq_algebraic_norms}. It applies to $J-2\alpha t$, and it gives, for $k\geq3$,
\begin{align}\label{eq_general_weyl}
\begin{split}
	J={}&4|D\II|^2+3|A|^2-6|d|^2+4\Big|L^s+\frac12\mathcal{P}\Big|^2+2\Big|V-\frac34K\Big|^2+2|Y|^2+\frac{1}{m^2\alpha^2}|\mathcal{F}|^2\\
	&+\Big(\frac{14}{k-1}\lambda+2\alpha\Big)t+6|\mathring{\Ric}_F|^2+\Big(\frac{28}{k-2}-6\Big)\inner{\mathring{\Ric}_F}{\mathring{T}}-14W^F(\II,\II)\\
	&-\frac{17}{4}C+\frac72c+\Big(\frac{7}{k(k-1)}-\frac{3}{2k}\Big)t^2.
\end{split}
\end{align}
When $k=2$, \eqref{eq_general_weyl} holds with the terms containing $\mathring{\Ric}_F$, $\mathring{T}$ and $W^F$ omitted.

We shall also need the following identity for the second fundamental form, which is \cite[Lemma 4.7]{horfocal}. Set
\begin{gather*}
	O_{\II}=\sum_{i,j,\ell}\inner{R^\perp(e_i,e_j)\II_{j\ell}}{\II_{i\ell}}=\frac12\inner{K}{V}+C,\\
	\omega_{\II}=\sum_{i,j,\ell,s}R^F_{ij\ell s}\inner{\II_{is}}{\II_{j\ell}}=-\sum_\beta R^F(S_\beta,S_\beta),
\end{gather*}
where the second expression for $O_{\II}$ follows from Lemma \ref{lem_first_order}(iii).

\begin{lemma}\label{lem_simons}
On the closed manifold $F$, $\int_FI\dmu_F=0$, where
\[
	I=C_D-|d|^2-O_{\II}+\trg(\Ric_FT)+\omega_{\II}.
\]
For every $k\geq3$,
\begin{align}\label{eq_simons_weyl}
	\trg(\Ric_FT)+\omega_{\II}=-W^F(\II,\II)+\frac{k}{k-2}\inner{\mathring{\Ric}_F}{\mathring{T}}+\frac{k}{k-1}\lambda t+\frac{t^2}{2(k-1)}.
\end{align}
\end{lemma}

\begin{proof}
The identity $\int_FI\dmu_F=0$ is the integrated Weitzenb\"ock formula for bundle-valued one-forms \cite[Proposition A.7]{drth}, applied to $\II$, regarded as a one-form on $F$ with values in $T^*F\otimes NF$. Then \eqref{eq_simons_weyl} follows from the decomposition of the curvature tensor of $F$ and Lemma \ref{lem_hess_two}(ii). Mutatis mutandis, the proof of \cite[Lemma 4.7]{horfocal} applies.
\end{proof}

\section{Two- and three-dimensional focal sets}\label{sec_lowdim}

The main goal of this section is to prove that $F$ is totally geodesic when its dimension is two or three. The identities below are those of \cite[Theorems 5.1 and 5.2]{horfocal}, with the additional term $2\alpha t$.

\begin{theorem}\label{thm_two}
If $k=2$, then, for every $q\geq1$,
\begin{align}\label{eq_two}
\begin{split}
	0=\int_F\Big(&4|D\II|^2+(14\lambda+2\alpha)t+\frac{19}{8}t^2+\frac{31}{4}\Big(c-\frac{t^2}{2}\Big)\\
	&+4\Big|L^s+\frac12\mathcal{P}\Big|^2+2\Big|V-\frac34K\Big|^2+2|Y|^2+\frac{1}{m^2\alpha^2}|\mathcal{F}|^2\Big)\dmu_F,
\end{split}
\end{align}
and every term of the integrand is nonnegative. Hence $F$ is totally geodesic.
\end{theorem}

\begin{proof}
By Theorem \ref{thm_general}, it suffices to show that, for $k=2$, $J$ is the integrand of \eqref{eq_two}. For $J-2\alpha t$, this is the computation in the proof of \cite[Theorem 5.1]{horfocal}. It starts from \eqref{eq_general_weyl} for $k=2$, and it uses the equality case of \eqref{eq_d_bound}, the identities \eqref{eq_algebraic_norms} and the fact that $S_\beta S_\gamma+S_\gamma S_\beta=G_{\beta\gamma}\Id$ for trace free symmetric endomorphisms of a two-dimensional space. Note that $G$ is the Gram matrix of $S_1,\ldots,S_q$ in the space of trace free symmetric endomorphisms of $T_pF$, which is a two-dimensional Euclidean space. Thus $G\geq0$ has rank at most two and trace $t$, so that $c=|G|^2\geq\frac12t^2$. Therefore every term of the integrand is nonnegative and, since $14\lambda+2\alpha>0$, it then follows that $t=0$, as desired.
\end{proof}

\begin{theorem}\label{thm_three}
If $k=3$, then, for every $q\geq1$,
\begin{align}\label{eq_three}
\begin{split}
	0=\int_F\Big(&5|A|^2-2|d|^2+(\lambda+2\alpha)t+6\Big|\mathring{\Ric}_F+\frac56\mathring{T}\Big|^2+\frac{13}{36}\big(C+t^2\big)+\frac{19}{9}c\\
	&+4\Big|L^s+\frac12\mathcal{P}\Big|^2+2\Big|V-\frac14K\Big|^2+2|Y|^2+\frac{1}{m^2\alpha^2}|\mathcal{F}|^2\Big)\dmu_F,
\end{split}
\end{align}
and the integrand is a sum of nonnegative terms, since $5|A|^2-2|d|^2\geq3|A|^2$ by \eqref{eq_d_bound}. Hence $F$ is totally geodesic.
\end{theorem}

\begin{proof}
By Theorem \ref{thm_general} and Lemma \ref{lem_simons}, we have $\int_F(J-4I)\dmu_F=0$. For $k=3$, mutatis mutandis, the proof of \cite[Theorem 5.2]{horfocal} shows that $J-2\alpha t-4I$ is the integrand of \eqref{eq_three} with $\lambda t$ in place of $(\lambda+2\alpha)t$. It uses \eqref{eq_general_weyl}, \eqref{eq_simons_weyl} with $W^F=0$, \eqref{eq_algebraic_norms} and the Cayley--Hamilton identity $\trg S^4=\frac12(\trg S^2)^2$ for trace free symmetric $3\times3$ matrices. Since $5|A|^2-2|d|^2\geq3|A|^2$ by \eqref{eq_d_bound}, every term of the integrand of \eqref{eq_three} is nonnegative, and therefore $t=0$, since $\lambda+2\alpha>0$.
\end{proof}

It follows from Theorems \ref{thm_two} and \ref{thm_three} that $F=\operatorname{MAX}(u)$ is totally geodesic whenever $k\in\{2,3\}$, in every dimension $n$.

\section{Totally geodesic maximum sets}\label{sec_criterion}

This section is devoted to the proof of Theorem \ref{thm_criterion}. The corresponding theorem for shrinking solitons is due to Li, Li, Ou, Qu, Wu and Zhou \cite{lloqwz}, and, mutatis mutandis, we will follow their argument. We keep the notation of Section \ref{sec_prelim} for $F=\operatorname{MAX}(u)$, $F_\varepsilon$, $r$, $\nu$, $\Sigma_r$, $S$ and $\psi$. From Subsection \ref{subsec_defect} to Subsection \ref{subsec_higher}, we will suppose that $k\geq2$ and $q\geq2$.

\subsection{The projection defect}\label{subsec_defect}

Let $\pi\colon F_\varepsilon\to F$ be the normal projection. Note that $\pi$ is a smooth submersion, by Proposition \ref{prop_structure}(iii). Let us set
\[
	\mathcal{N}=\ker\mathrm{d}\pi,\qquad \mathcal{H}=\mathcal{N}^\perp,
\]
and let us denote by $\Pi_{\mathcal{N}}$ and $\Pi_{\mathcal{H}}$ the orthogonal projections onto these two distributions, which are smooth on $F_\varepsilon$ and which, at the points of $F$, are $NF$ and $TF$. For each $x\in F$, we will denote by $\Sph_x$ the unit sphere of $N_xF$ and, for $0<r<\varepsilon$, by $\varphi_r\colon\Sph_x\to M$ the map $\varphi_r(\theta)=\exp_x(r\theta)$, whose image we will call a normal sphere. The normal spheres are tangent to a distribution $\mathcal{V}$ of rank $q-1$ on $F_\varepsilon\setminus F$, and there we have $\mathcal{N}=\R\nu\oplus\mathcal{V}$ and $T\Sigma_r=\mathcal{H}\oplus\mathcal{V}$. The orthogonal projection onto $\mathcal{V}$ will be denoted by $\Pi_{\mathcal{V}}$. With this notation, we define
\begin{align}\label{eq_defect}
	E=\frac{1}{2m\alpha}P-\frac12\Pi_{\mathcal{H}},\qquad v=\trg_{\mathcal{V}}E,\qquad Z=\Pi_{\mathcal{H}}E|_{\mathcal{V}}.
\end{align}
Note that $E$ is a smooth symmetric tensor on $F_\varepsilon$ which, by Proposition \ref{prop_structure}, vanishes on $F$ and satisfies $E\nu=0$. Recalling the definition of $E$, we see that \eqref{eq_shape} reads
\begin{align}\label{eq_shape_E}
	S=\frac{\psi'}{\psi}\big(\Id-\Pi_{\mathcal{H}}-2E\big)\qquad\text{on } T\Sigma_r.
\end{align}

\begin{lemma}\label{lem_square}
On $F_\varepsilon\setminus F$, the function $v$ is the squared norm of $E$, that is, $|E|^2=v$.
\end{lemma}

\begin{proof}
By \eqref{eq_constants}, we have $|P|^2=km^2\alpha^2$ and $\trg P=km\alpha$. Since $\Pi_{\mathcal{H}}$ is an orthogonal projection of rank $k$, we also have $|\Pi_{\mathcal{H}}|^2=k$, and therefore
\[
	|E|^2=\frac k4-\frac{1}{2m\alpha}\trg_{\mathcal{H}}P+\frac k4=\frac{1}{2m\alpha}\big(\trg P-\trg_{\mathcal{H}}P\big)=\frac{1}{2m\alpha}\trg_{\mathcal{V}}P,
\]
where we used $P(\nu,\cdot)=0$ in the last equality. But $\Pi_{\mathcal{H}}$ vanishes on $\mathcal{V}$, and therefore $v=\frac{1}{2m\alpha}\trg_{\mathcal{V}}P$, as desired.
\end{proof}

On each sphere $\Sph_x$, we consider the metric $h_r=\psi(r)^{-2}\varphi_r^*g$. Let $g_0$ be the round metric of $\Sph_x$, $|\Sph^{q-1}|$ its volume and $\mathrm{d}\theta$ the round probability measure. We define a function $\chi>0$ on $F_\varepsilon\setminus F$ by means of $\mathrm{d}V_{h_r}=|\Sph^{q-1}|\,\chi(\exp_x(r\theta))\,\mathrm{d}\theta$. Finally, given a vector $X$ tangent to $F$ at $x$, its horizontal lift is denoted by $\hat X$, namely $\hat X$ is the section of $\mathcal{H}$ along $\pi^{-1}(x)$ which $\mathrm{d}\pi$ maps to $X$.

\begin{lemma}\label{lem_evolution}
On $F_\varepsilon\setminus F$:
\begin{enumerate}[label=\textup{(\roman*)},leftmargin=*]
\item $\partial_rh_r=-4\frac{\psi'}{\psi^3}\varphi_r^*\big(g(E\cdot,\cdot)\big)$ and $\partial_r\log\chi=-2\frac{\psi'}{\psi}v$, and $h_r\to g_0$ and $\chi\to1$ as $r\to0$;
\item $\partial_r\,g(\hat X,\hat Y)=-4\frac{\psi'}{\psi}g(E\hat X,\hat Y)$ for $X,Y$ tangent to $F$ at $x$, and the volume density of the metric $(X,Y)\mapsto g(\hat X,\hat Y)$ on $T_xF$ with respect to $g_F$ is $1/\chi$.
\end{enumerate}
\end{lemma}

\begin{proof}
(i) Let $\zeta$ be tangent to $\Sph_x$ at $\theta$. Then $w=\mathrm{d}\varphi_r(\zeta)$ is the Jacobi field along the geodesic $r\mapsto\exp_x(r\theta)$ which vanishes at $r=0$ and has derivative $\zeta$ there. Note that $w$ is tangent to $\mathcal{V}$ and that $[\nu,w]=0$, whence $\nabla_\nu w=Sw$. Given two such fields $w_1$ and $w_2$, it follows from \eqref{eq_shape_E} that $\partial_rg(w_1,w_2)=2g(Sw_1,w_2)=2\frac{\psi'}{\psi}\big(g(w_1,w_2)-2g(Ew_1,w_2)\big)$. Consequently, $\partial_r\big(\psi^{-2}g(w_1,w_2)\big)=-4\frac{\psi'}{\psi^3}g(Ew_1,w_2)$, which is the first identity. Taking the trace of this identity, we also have $\partial_r\log\mathrm{d}V_{h_r}=\frac12\trg_{h_r}(\partial_rh_r)=-2\frac{\psi'}{\psi}v$. As for the limits, they follow from $w=r\zeta+O(r^3)$, in a parallel frame along the geodesic, and from $\psi(r)=r+O(r^3)$, which give $h_r\to g_0$ and $\chi\to1$ as $r\to0$.

(ii) We may suppose that $X$ and $Y$ are vector fields on $F$. The field $\hat X$ is projectable and $\nu$ is tangent to the fibers of $\pi$, so $[\nu,\hat X]$ is a section of $\mathcal{N}$ and $g(\nabla_\nu\hat X,\hat Y)=g(\nabla_{\hat X}\nu,\hat Y)=g(S\hat X,\hat Y)$. The first assertion is then an immediate consequence of \eqref{eq_shape_E}. For the second one, let us compute in two ways the density of the volume measure of $\Sigma_r$ with respect to the product of the measure of $F$ and of $\mathrm{d}\theta$. On the one hand, $\mathcal{H}$ and $\mathcal{V}$ are orthogonal, and therefore this density is the product of the volume density of the metric $(X,Y)\mapsto g(\hat X,\hat Y)$ and of the density $|\Sph^{q-1}|\psi^{q-1}\chi$ of the measure of the normal spheres with respect to $\mathrm{d}\theta$. On the other hand, it is $|\Sph^{q-1}|\psi^{q-1}$, by Proposition \ref{prop_structure}(iv). Hence the volume density of the metric $(X,Y)\mapsto g(\hat X,\hat Y)$ is $1/\chi$, as claimed.
\end{proof}

From now on, for $x\in F$ and $\theta\in\Sph_x$, the tangent spaces along the geodesic $r\mapsto\exp_x(r\theta)$ will be identified with $T_xM=T_xF\oplus N_xF$ by means of parallel transport, and $S_\theta$ will denote the shape operator of $F$ in the direction of $\theta$. Note that all the estimates $O(r^j)$ below are uniform on the unit normal bundle of $F$, since it is compact.

\begin{lemma}\label{lem_first_jet}
In the splitting $T_xM=T_xF\oplus N_xF$, we have
\[
	E(\exp_x(r\theta))=\frac r2\begin{pmatrix}S_\theta&0\\0&0\end{pmatrix}+O(r^2),\qquad \Pi_{\mathcal{H}}(\exp_x(r\theta))=\Pi_T+O(r^2).
\]
A fortiori, $E=O(r^2)$ whenever $F$ is totally geodesic.
\end{lemma}

\begin{proof}
Recall that, by Lemma \ref{lem_first_order}(i), $\nabla_\theta P$ is $m\alpha S_\theta$ on $T_xF$ and has no other nonzero component. Let us fix $\zeta\in N_xF$ orthogonal to $\theta$, and let $w$ be the corresponding Jacobi field of the proof of Lemma \ref{lem_evolution}, so that $w=r\zeta+O(r^3)$. Such fields, together with $\nu$, span $\mathcal{N}$ along the geodesic, whence $\Pi_{\mathcal{N}}=\Pi_N+O(r^2)$ and $\Pi_{\mathcal{H}}=\Pi_T+O(r^2)$. Since $E$ vanishes at $x$, the expansion of $E$ then follows from Taylor's formula.
\end{proof}

\subsection{A scalar balance}\label{subsec_balance}

In this subsection, we relate $E$ to the scalar curvature $R_{h_r}$ of the metric $h_r$. If $q=2$, we set $R_{h_r}=0$.

\begin{lemma}\label{lem_fiber}
On $F_\varepsilon\setminus F$,
\begin{align}\label{eq_fiber}
	\psi^2\trg_{\mathcal{V}}\Ric_\Sigma=(q-1)(q-2)+\big(2\lambda\psi^2-2(q-2)\big)v-2\psi\psi'\,\partial_rv-8(\psi')^2|Z|^2.
\end{align}
\end{lemma}

\begin{proof}
Let $w$ be a Jacobi field as in the proof of Lemma \ref{lem_evolution}. By \eqref{eq_shape_E}, we have $\nabla_\nu w=\frac{\psi'}{\psi}(w-2Ew)$. Since $E\nu=0$, the component of $\nabla_\nu w$ orthogonal to $\mathcal{V}$ is $-2\frac{\psi'}{\psi}Zw$. Consequently, with respect to $\mathcal{H}\oplus\mathcal{V}$,
\[
	\nabla_\nu\Pi_{\mathcal{V}}=-2\frac{\psi'}{\psi}\begin{pmatrix}0&Z\\Z^*&0\end{pmatrix},\qquad \partial_rv=\trg_{\mathcal{V}}(\nabla_\nu E)-4\frac{\psi'}{\psi}|Z|^2.
\]
By \eqref{eq_defect}, $P=m\alpha\Pi_{\mathcal{H}}+2m\alpha E$. Since $\nabla_\nu\Pi_{\mathcal{H}}=-\nabla_\nu\Pi_{\mathcal{V}}$ has no component on $\mathcal{V}\times\mathcal{V}$, it follows that
\[
	\trg_{\mathcal{V}}P=2m\alpha v,\qquad \trg_{\mathcal{V}}(\nabla_\nu P)=2m\alpha\Big(\partial_rv+4\frac{\psi'}{\psi}|Z|^2\Big),
\]
and \eqref{eq_fiber} follows by taking the trace of \eqref{eq_ric_sigma} over $\mathcal{V}$ and using these two identities.
\end{proof}

We denote by $\mathrm{d}\mu_0$ the product of the measure of $F$ and of $\mathrm{d}\theta$, divided by $|F|$. This is a probability measure on the unit normal bundle of $F$. For the remainder of this section, all integrals are taken over the unit normal bundle with respect to $\mathrm{d}\mu_0$, for a fixed value of $r$. Since the scalar curvature of the normal sphere $\varphi_r(\Sph_x)$ is $\psi^{-2}R_{h_r}$, we set
\[
	\Gamma(r)=\int\big(\psi^{-2}R_{h_r}-\trg_{\mathcal{V}}\Ric_\Sigma\big)\chi.
\]
By Lemma \ref{lem_evolution}, $\partial_r\chi=-2\frac{\psi'}{\psi}v\chi$, which implies that $\int(\partial_rv)\chi=\frac{\mathrm{d}}{\mathrm{d}r}\int v\chi+2\frac{\psi'}{\psi}\int v^2\chi$. Therefore, multiplying \eqref{eq_fiber} by $\chi$, integrating and using the definition of $\Gamma(r)$, we have
\begin{align}\label{eq_balance}
\begin{split}
	\int\big(R_{h_r}-(q-1)(q-2)\big)\chi={}&\big(2\lambda\psi^2-2(q-2)\big)\int v\chi-2\psi\psi'\frac{\mathrm{d}}{\mathrm{d}r}\int v\chi\\
	&-4(\psi')^2\int v^2\chi-8(\psi')^2\int|Z|^2\chi+\psi^2\,\Gamma(r).
\end{split}
\end{align}
For solitons, the corresponding identity is \eqref{eq_balance} with $\psi(r)=r$; see \cite[Section 6]{lloqwz}.

\subsection{A spherical estimate}\label{subsec_sphere}

Let $h$ be a metric on the sphere $\Sph^{q-1}$, $q\geq3$, and consider the functional
\[
	\Phi(h)=\Big(\frac{1}{|\Sph^{q-1}|}\int_{\Sph^{q-1}}\mathrm{d}V_h\Big)^{-\frac{q-3}{q-1}}\frac{1}{|\Sph^{q-1}|}\int_{\Sph^{q-1}}R_h\,\mathrm{d}V_h,
\]
so that $\Phi(g_0)=(q-1)(q-2)$ for the round metric $g_0$ of curvature one. If $q=2$, we set $\Phi=0$. The metric $g_0$ is a critical point of $\Phi$, and $D^2\Phi$ denotes the second variation of $\Phi$ at $g_0$. The following estimate is \cite[Lemma 7.1]{lloqwz}. For the convenience of the reader, we include a proof.

\begin{lemma}\label{lem_sphere}
Let $N$ be a trace free symmetric $2$-tensor on the unit sphere $\Sph^{q-1}\subset\R^q$ such that, after extending $N$ by zero in the radial direction, its components in the canonical basis of $\R^q$ are restrictions of polynomials of degree at most $d$. Then
\[
	D^2\Phi(N,N)\geq-\frac12d(d+q-2)\int_{\Sph^{q-1}}|N|^2\,\mathrm{d}\theta.
\]
\end{lemma}

\begin{proof}
If $q=2$, there is nothing to prove, and, if $q=3$, then $\Phi$ is constant by the Gauss--Bonnet theorem, so that the inequality is trivially satisfied. We will therefore assume that $q\geq4$. The first variation of $\Phi$ at a metric $h$ in the direction $N$ is the integral of the inner product of $N$ with a multiple of $-\Ric_h+\frac12R_hh-\frac{q-3}{2(q-1)}\frac{\int R_h\,\mathrm{d}V_h}{\int\mathrm{d}V_h}\,h$. This tensor vanishes at $g_0$. We differentiate it at $g_0$ in the direction $N$. Since $N$ is trace free, the terms which are multiples of $g_0$ have zero inner product with $N$, and we are left with
\[
	D^2\Phi(N,N)=-\int_{\Sph^{q-1}}\inner{\Ric'(N)-(q-2)N}{N}\,\mathrm{d}\theta,
\]
where $\Ric'(N)$ denotes the derivative of the Ricci tensor at $g_0$ in the direction $N$. Let $\nabla$ and $\divg$ be taken with respect to $g_0$. For trace free $N$ on the unit sphere, we have $\Ric'(N)=\frac12\big(\nabla^*\nabla N+2(q-1)N\big)+\frac12\mathcal{L}_{(\divg N)^\sharp}g_0$; see \cite[Theorem 1.174]{besse}. Integrating by parts, we get
\[
	D^2\Phi(N,N)=-\frac12\int_{\Sph^{q-1}}\big(|\nabla N|^2+2|N|^2\big)\,\mathrm{d}\theta+\int_{\Sph^{q-1}}|\divg N|^2\,\mathrm{d}\theta\geq-\frac12\int_{\Sph^{q-1}}\big(|\nabla N|^2+2|N|^2\big)\,\mathrm{d}\theta.
\]
Now we regard $N$ as a symmetric bilinear form on $\R^q$ such that $N(\theta,\cdot)=0$ at each $\theta\in\Sph^{q-1}$, and we denote by $N_{\beta\gamma}$ its components in the canonical basis. Differentiating $N(\theta,\cdot)=0$ along the sphere, it is easy to see that the sum over $\beta,\gamma$ of the squared norms of the gradients of the functions $N_{\beta\gamma}$ on $\Sph^{q-1}$ equals $|\nabla N|^2+2|N|^2$. Moreover, each $N_{\beta\gamma}$ is the restriction to $\Sph^{q-1}$ of a polynomial of degree at most $d$. It is therefore a sum of spherical harmonics of degree at most $d$, and the corresponding eigenvalues of the Laplacian do not exceed $d(d+q-2)$. Consequently,
\[
	\int_{\Sph^{q-1}}\big(|\nabla N|^2+2|N|^2\big)\,\mathrm{d}\theta\leq d(d+q-2)\int_{\Sph^{q-1}}|N|^2\,\mathrm{d}\theta,
\]
as desired.
\end{proof}

\subsection{The higher normal derivatives}\label{subsec_higher}

In this subsection, we prove the following proposition, which is the main step in the proof of Theorem \ref{thm_criterion}.

\begin{proposition}\label{prop_higher}
If $q\geq2$ and $E=O(r^2)$, then $E=O(r^j)$ for every $j\geq1$.
\end{proposition}

\begin{proof}
Suppose by contradiction that the assertion is false. Then there is a largest integer $d\geq2$ such that $E=O(r^d)$. In orthonormal frames of $TF$ and $NF$ extended by parallel transport along the normal geodesics of $F$, the components of $E$ are smooth functions on $F_\varepsilon$ which vanish to order $d-1$ along $F$. It follows that
\begin{align}\label{eq_leading}
	E(\exp_x(r\theta))=r^dE_d(x,\theta)+O(r^{d+1}),\qquad E_d=\begin{pmatrix}H_d&Z_d\\Z_d^*&N_d\end{pmatrix},
\end{align}
where $E_d(x,\cdot)$ is a homogeneous polynomial of degree $d$ on $N_xF$ with values in the symmetric endomorphisms of $T_xM=T_xF\oplus N_xF$. By the choice of $d$, $E_d$ does not vanish identically on the unit normal bundle. Moreover, the expansion can be differentiated in $x$ and $\theta$. Note that $E\nu=0$ implies that $N_d\theta=0$ and $Z_d\theta=0$. By Lemma \ref{lem_first_jet}, $H_d$, $Z_d$ and $N_d$ are also the leading coefficients of the blocks of $E$ with respect to $\mathcal{H}\oplus\mathcal{N}$. In particular,
\begin{align}\label{eq_Ed_norm}
	\int|E_d|^2=\int|H_d|^2+\int|N_d|^2+2\int|Z_d|^2>0.
\end{align}
It follows from Lemma \ref{lem_square} and \eqref{eq_leading} that $v=r^{2d}|E_d|^2+O(r^{2d+1})$ and $|Z|^2=r^{2d}|Z_d|^2+O(r^{2d+1})$. Since $\psi'/\psi=r^{-1}+O(r)$, we can integrate the identities of Lemma \ref{lem_evolution}, which yields
\begin{align}\label{eq_expansions}
\begin{split}
	\chi&=1-\frac{r^{2d}}{d}|E_d|^2+O(r^{2d+1}),\qquad h_r=g_0-\frac4dr^dN_d+O(r^{d+1}),\\
	g(\hat X,\hat Y)&=g_F(X,Y)-\frac4dr^d\inner{H_dX}{Y}+O(r^{d+1}),
\end{split}
\end{align}
where $N_d$ is restricted to the tangent spaces of $\Sph_x$. Here we used that $\chi$ and $|Z|^2=|\Pi_{\mathcal{H}}E\,\Pi_{\mathcal{N}}|^2$ are smooth functions on $F_\varepsilon$. Indeed, $\mathcal{H}$ and $\mathcal{N}$ are smooth distributions, and $1/\chi$ is the volume density of Lemma \ref{lem_evolution}(ii). Furthermore, the components of $\psi^2h_r=\varphi_r^*g$ in the normal coordinates $r\theta$ are smooth, and therefore the expansion of $h_r$ holds together with its derivatives along $\Sph_x$, uniformly in $x$. On the other hand, $v=\trg_{\mathcal{V}}E=r^d\,\trg N_d+O(r^{d+1})$. Comparing this with the expansion of $v$ above, we conclude that $\trg N_d=0$. Finally, since $N_d\theta=0$, the extension of the restriction of $N_d$ by zero in the radial direction is $N_d$ itself, whose components are polynomials of degree $d$.

First, we will estimate $\Gamma(r)$ from below. Since $g_0$ is a critical point of $\Phi$, it follows from \eqref{eq_expansions} that $\Phi(h_r)=(q-1)(q-2)+\frac{8}{d^2}r^{2d}D^2\Phi(N_d,N_d)+O(r^{2d+1})$ on each $\Sph_x$. Moreover, by the definition of $\Phi$, the integral $\int_{\Sph_x}R_{h_r}\chi\,\mathrm{d}\theta$ is the product of $\Phi(h_r)$ and of the power $\frac{q-3}{q-1}$ of $\int_{\Sph_x}\chi\,\mathrm{d}\theta$. We expand this product using \eqref{eq_expansions}, apply Lemma \ref{lem_sphere} on each $\Sph_x$ and average over $F$. The result is
\[
	\int\big(R_{h_r}-(q-1)(q-2)\big)\chi\geq r^{2d}\Big(\frac{2(q-2)}{d}\int|E_d|^2-\frac{4(d+q-2)}{d}\int|N_d|^2\Big)+O(r^{2d+1}),
\]
and this inequality also holds when $q=2$, because in this case the left-hand side vanishes and $N_d=0$. We also have
\begin{gather*}
	\int v\chi=r^{2d}\int|E_d|^2+O(r^{2d+1}),\qquad \frac{\mathrm{d}}{\mathrm{d}r}\int v\chi=2d\,r^{2d-1}\int|E_d|^2+O(r^{2d}),\\
	\int|Z|^2\chi=r^{2d}\int|Z_d|^2+O(r^{2d+1}),\qquad \int v^2\chi=O(r^{4d}).
\end{gather*}
Since $\psi^2=r^2+O(r^4)$, $\psi\psi'=r+O(r^3)$ and $(\psi')^2=1+O(r^2)$, the term of \eqref{eq_balance} which contains $\lambda\psi^2$ is $O(r^{2d+2})$. Inserting the expansions above into \eqref{eq_balance}, we conclude that
\begin{align}\label{eq_gamma_lower}
	\liminf_{r\to0}r^{-(2d-2)}\Gamma(r)\geq\Big(4d+2(q-2)+\frac{2(q-2)}{d}\Big)\int|E_d|^2+8\int|Z_d|^2-\frac{4(d+q-2)}{d}\int|N_d|^2.
\end{align}

Next, we will estimate $\Gamma(r)$ from above. We work inside $\Sigma_r$, whose Levi-Civita connection is denoted by $\nabla^\Sigma$. Let
\[
	B_{\mathcal{H}}(X,Y)=\frac12\big(\nabla^\Sigma_XY+\nabla^\Sigma_YX\big)^{\mathcal{V}}
\]
for sections $X,Y$ of $\mathcal{H}$, and define $B_{\mathcal{V}}$ in the same way, with the roles of $\mathcal{H}$ and $\mathcal{V}$ exchanged. Consider the sum of the sectional curvatures of $\Sigma_r$ on the planes spanned by a vector of an orthonormal basis of $\mathcal{V}$ and a vector of an orthonormal basis $\{X_i\}$ of $\mathcal{H}$. Since $\mathcal{V}$ is integrable, the formula of Walczak \cite{walczak} for the orthogonal distributions $\mathcal{H}$ and $\mathcal{V}$ on $\Sigma_r$ shows that this sum equals
\[
	\divg_\Sigma\big(\trg B_{\mathcal{H}}+\trg B_{\mathcal{V}}\big)-|B_{\mathcal{H}}|^2-|B_{\mathcal{V}}|^2+|\trg B_{\mathcal{H}}|^2+|\trg B_{\mathcal{V}}|^2+\frac14\sum_{i,j}\big|[X_i,X_j]^{\mathcal{V}}\big|^2.
\]
On the other hand, by the Gauss equation of the normal sphere in $\Sigma_r$, the scalar curvature of the normal sphere is $\trg_{\mathcal{V}}\Ric_\Sigma$ minus this sum plus $|\trg B_{\mathcal{V}}|^2-|B_{\mathcal{V}}|^2$. Combining these two facts, we have
\begin{align}\label{eq_walczak}
	\psi^{-2}R_{h_r}-\trg_{\mathcal{V}}\Ric_\Sigma=|B_{\mathcal{H}}|^2-\frac14\sum_{i,j}\big|[X_i,X_j]^{\mathcal{V}}\big|^2-|\trg B_{\mathcal{H}}|^2-\divg_\Sigma\big(\trg B_{\mathcal{H}}+\trg B_{\mathcal{V}}\big).
\end{align}
We now compute $\trg B_{\mathcal{H}}$. Let $w$ be a section of $\mathcal{V}$, and suppose that $X_1,\ldots,X_k$ are horizontal lifts of vector fields on $F$ which are orthonormal at a given point. Since $[w,X_i]$ is a section of $\mathcal{V}$, at this point the derivative in the direction of $w$ of the logarithm of the volume density of Lemma \ref{lem_evolution}(ii) is $\sum_ig(\nabla_wX_i,X_i)=\sum_ig(\nabla_{X_i}w,X_i)=-g(w,\trg B_{\mathcal{H}})$. Since this density is $1/\chi$, it follows that $\trg B_{\mathcal{H}}$ is the component of $\nabla\log\chi$ in $\mathcal{V}$. Moreover, by Proposition \ref{prop_structure}(iv), the volume measure of $\Sigma_r$ is a constant multiple of $\mathrm{d}\mu_0$. Thus, for every vector field tangent to $\Sigma_r$, the integral of $\chi$ times its divergence is minus the integral of $\chi$ times its inner product with $\nabla\log\chi$. We multiply \eqref{eq_walczak} by $\chi$ and integrate. The terms containing $|\trg B_{\mathcal{H}}|^2$ cancel, and we arrive at
\begin{align}\label{eq_gamma_walczak}
	\Gamma(r)=\int\chi\Big(|B_{\mathcal{H}}|^2-\frac14\sum_{i,j}\big|[X_i,X_j]^{\mathcal{V}}\big|^2+g\big(\nabla\log\chi,\trg B_{\mathcal{V}}\big)\Big).
\end{align}
It remains to estimate each term of \eqref{eq_gamma_walczak}. Given vector fields $X,Y$ on $F$ and a section $w$ of $\mathcal{V}$, the Koszul formula yields
\[
	2g\big(B_{\mathcal{H}}(\hat X,\hat Y),w\big)=-w\big(g(\hat X,\hat Y)\big),
\]
because $[w,\hat X]$ and $[w,\hat Y]$ are sections of $\mathcal{V}$. A unit vector of $\mathcal{V}$ has the form $r^{-1}(1+O(r^2))\,\mathrm{d}\varphi_r(\zeta)$, where $\zeta$ is a $g_0$-unit vector tangent to $\Sph_x$. Thus \eqref{eq_expansions} implies that
\[
	|B_{\mathcal{H}}|^2=\frac{4}{d^2}r^{2d-2}|\nabla^{\Sph}H_d|^2+O(r^{2d-1}),
\]
where $\nabla^{\Sph}H_d$ denotes the gradient along $\Sph_x$ of the components of $H_d(x,\cdot)$ in an orthonormal basis of $T_xF$. Now consider the last term of \eqref{eq_gamma_walczak}. By \eqref{eq_expansions}, $\log\chi=O(r^{2d})$. Since the horizontal lifts have bounded components in $x$ and $\theta$, the component of $\nabla\log\chi$ in $\mathcal{H}$ is $O(r^{2d})$ as well. On the other hand, $\trg B_{\mathcal{V}}$ is the trace over $\mathcal{V}$ of the component in $\mathcal{H}$ of the second fundamental form in $M$ of the image under $\exp_x$ of the ball of radius $\varepsilon$ of $N_xF$. The geodesics of $M$ through $x$ which are normal to $F$ are contained in this submanifold. Hence its second fundamental form vanishes at $x$, and we have $\trg B_{\mathcal{V}}=O(r)$. Since the second term of \eqref{eq_gamma_walczak} is nonpositive, we can discard it, and we conclude that
\begin{align}\label{eq_gamma_upper}
	\limsup_{r\to0}r^{-(2d-2)}\Gamma(r)\leq\frac{4}{d^2}\int|\nabla^{\Sph}H_d|^2\leq\frac{4(d+q-2)}{d}\int|H_d|^2,
\end{align}
where the second inequality holds because the components of $H_d(x,\cdot)$ are restrictions to $\Sph_x$ of polynomials of degree $d$, and such restrictions are sums of spherical harmonics with eigenvalues at most $d(d+q-2)$.

Finally, we compare \eqref{eq_gamma_lower} with \eqref{eq_gamma_upper}. Using \eqref{eq_Ed_norm}, we get
\[
	\Big(4(d-1)+2(q-2)\Big(1-\frac1d\Big)\Big)\Big(\int|H_d|^2+\int|N_d|^2\Big)+\Big(8d+8+4(q-2)\Big(1+\frac1d\Big)\Big)\int|Z_d|^2\leq0.
\]
Since $d\geq2$, both coefficients are positive, which is evidently absurd, in view of \eqref{eq_Ed_norm}.
\end{proof}

\begin{remark}\label{rem_first_mode}
Note that, for $d=1$, the coefficient of $\int|H_d|^2+\int|N_d|^2$ in the last inequality vanishes. On the other hand, Lemma \ref{lem_first_jet} shows that the coefficient of degree one of $E$ is given by the second fundamental form of $F$. This is the reason why we assume that $F$ is totally geodesic in Theorem \ref{thm_criterion}; see also \cite[Remark 8.1]{lloqwz}.
\end{remark}

\subsection{Proof of Theorem \ref{thm_criterion}}\label{subsec_proof_criterion}

We start by showing that rigidity follows from the identities $P^2=m\alpha P$ and $\nabla_{\nabla u}P=0$.

\begin{proposition}\label{prop_splitting}
Suppose that $P^2=m\alpha P$ and $\nabla_{\nabla u}P=0$ on a nonempty open subset of the interior of $M$. Then $M$ is rigid, and its universal cover is $\bar N^k\times\Sph^q_+(\alpha)$, where $\bar N^k$ is a compact simply connected Einstein manifold with $\Ric_{\bar N}=\lambda g_{\bar N}$.
\end{proposition}

\begin{proof}
It is well known that $g$ and $u$ are real analytic in harmonic coordinates in the interior of $M$; see \cite[Proposition 2.4]{hpw1}. Consequently, the tensors $P^2-m\alpha P$ and $\nabla_{\nabla u}P$ are real analytic as well. Since the interior of $M$ is connected, these tensors vanish on it, and by continuity they vanish on $M$. Now we differentiate $\nabla^2u=\frac um(P-m\alpha g)$, commute derivatives by means of \eqref{eq_ricci_identity} and contract with $\nabla u$. Using $P(\nabla u,\cdot)=0$, we get
\begin{align}\label{eq_radial}
\begin{split}
	\Rm(\nabla u,X,\nabla u,Y)={}&\frac{|\nabla u|^2}{m}(m\alpha g-P)(X,Y)-\alpha\,\mathrm{d}u(X)\mathrm{d}u(Y)\\
	&-\frac um(\nabla_{\nabla u}P)(X,Y)-\Big(\frac um\Big)^2(P^2-m\alpha P)(X,Y)
\end{split}
\end{align}
for all $X,Y$. The last two terms vanish, and $m\alpha g-P=\lambda g-\Ric$. Therefore
\[
	\Rm(\nabla u,X,\nabla u,Y)=\frac{|\nabla u|^2}{m}\big(\lambda g-\Ric\big)(X,Y)-\alpha\,\mathrm{d}u(X)\mathrm{d}u(Y)\qquad\text{for all } X,Y,
\]
which says that the tensor $Q$ defined in \cite[Section 3]{hpw2} is radially flat. This is the hypothesis of \cite[Theorem 1.11]{hpw2}. For $X=Y$ orthogonal to $\nabla u$, it is the identity (1.12) there. We apply this theorem to the universal cover of $M$, endowed with the lifted metric and potential function. The universal cover is a complete simply connected nontrivial $m$-quasi-Einstein manifold with boundary, constant scalar curvature, $\lambda>0$ and bounded curvature. It then follows that $M$ is rigid, and it remains to identify the factors. By \cite[Proposition 2.6]{hpw2}, the universal cover of $M$ is the product of an Einstein manifold $\bar N$ with $\Ric_{\bar N}=\lambda g_{\bar N}$ and of one of the Einstein examples with boundary and $\lambda>0$ of \cite[Table 2.1]{hpw2}, namely a hemisphere or an interval, on which the lifted potential function is defined. Here $\bar N$ is a point if $M$ is Einstein. Since $\Ric=\lambda\Pi_T+\rho\Pi_N$ at the points of $F$ and $\rho\neq\lambda$, $\bar N$ has dimension $k$, while the second factor has dimension $q$ and Ricci curvature $\rho=(q-1)\alpha$. By \eqref{eq_u_cos}, the second factor is $\Sph^q_+(\alpha)$ if $q\geq2$, and an interval of length $\pi/\sqrt{\alpha}$ if $q=1$. Finally, $\bar N$ is complete and simply connected, and it is compact by the theorem of Myers.
\end{proof}

\begin{proof}[Proof of Theorem \ref{thm_criterion}]
Suppose that $F=\operatorname{MAX}(u)$ is totally geodesic. If $k=0$, then $P=0$ by \eqref{eq_constants}, and the result follows from Proposition \ref{prop_splitting}. Since $k\neq1$, from now on we will assume that $k\geq2$.

Suppose first that $q=1$. In this case, \eqref{eq_shape} gives $S=0$ on $F_\varepsilon\setminus F$, that is, $P=m\alpha(g-\mathrm{d}r\otimes\mathrm{d}r)$. Since $\nabla_\nu\nu=0$, we have $P^2=m\alpha P$ and $\nabla_{\nabla u}P=0$ on $F_\varepsilon\setminus F$. Thus $M$ is rigid by Proposition \ref{prop_splitting}.

Now consider the case $q\geq2$. Since $F$ is totally geodesic, Lemma \ref{lem_first_jet} and Proposition \ref{prop_higher} imply that $E=O(r^j)$ for every $j\geq1$. Using $P=m\alpha\Pi_{\mathcal{H}}+2m\alpha E$ and $\Pi_{\mathcal{H}}^2=\Pi_{\mathcal{H}}$, we compute
\[
	P^2-m\alpha P=2m^2\alpha^2\big(\Pi_{\mathcal{H}}E+E\Pi_{\mathcal{H}}-E+2E^2\big),
\]
and therefore $P^2-m\alpha P=O(r^j)$ for every $j$. Next, we consider $\nabla_{\nabla u}P$. By \eqref{eq_u_cos}, we have $\nabla u=-\alpha U\psi\,\nu$ on $F_\varepsilon\setminus F$. Here $\psi\nu$ is the product of the radial vector field $r\nu$, which is smooth on $F_\varepsilon$, and of the smooth even function $\psi(r)/r$. Moreover, the proof of Lemma \ref{lem_fiber} shows that
\[
	\nabla_{\psi\nu}P=m\alpha\nabla_{\psi\nu}\Pi_{\mathcal{H}}+2m\alpha\nabla_{\psi\nu}E,\qquad \nabla_{\psi\nu}\Pi_{\mathcal{H}}=-\nabla_{\psi\nu}\Pi_{\mathcal{V}}=2\psi'\begin{pmatrix}0&Z\\Z^*&0\end{pmatrix},
\]
where $Z=\Pi_{\mathcal{H}}E|_{\mathcal{V}}$. Analogously, $\nabla_{\nabla u}P=O(r^j)$ for every $j$. Consequently, all the derivatives of the smooth tensors $P^2-m\alpha P$ and $\nabla_{\nabla u}P$ vanish at the points of $F$. Since these tensors are real analytic in the interior of $M$ \cite[Proposition 2.4]{hpw1}, they vanish on a neighborhood of $F$. Therefore $M$ is rigid by Proposition \ref{prop_splitting}.
\end{proof}

\section{Proofs of Theorems \ref{cor_alldim} and \ref{thm_five}}\label{sec_level}

We are now ready to prove Theorems \ref{cor_alldim} and \ref{thm_five}.

\begin{proof}[Proof of Theorem \ref{cor_alldim}]
By Theorems \ref{thm_two} and \ref{thm_three}, $F=\operatorname{MAX}(u)$ is totally geodesic. Theorem \ref{thm_criterion} then tells us that $M$ is rigid, with universal cover $\bar N^k\times\Sph^{n-k}_+(\alpha)$, where $\bar N^k$ is a compact simply connected Einstein manifold with $\Ric_{\bar N}=\lambda g_{\bar N}$. In dimension $k\in\{2,3\}$ such a manifold has constant sectional curvature, so $\bar N^k$ is the sphere $\Sph^k(\lambda/(k-1))$. It remains to describe the deck transformations. They are isometries, and they preserve the lifted potential function and the eigenspaces of the Ricci tensor, which are tangent to the two factors because $\rho\neq\lambda$. Consequently, each deck transformation is the product of an isometry of $\Sph^k(\lambda/(k-1))$ and an isometry of the hemisphere. The latter fixes the pole of the hemisphere, because the lifted potential function attains its maximum there. Since the group of deck transformations acts freely, it follows that it acts freely on $\Sph^k(\lambda/(k-1))$. This group is finite, because the universal cover is compact. Now suppose that $k=2$. Every isometry of $\Sph^2(\lambda)$ which preserves the orientation has a fixed point. Thus a deck transformation which preserves the orientation of $\Sph^2(\lambda)$ is the identity. But the product of two deck transformations which reverse the orientation preserves it, whence the group has order at most two, as claimed. Finally, suppose that $n=6$. By \eqref{eq_quantization}, we have $k\in\{0,2,3,4,5\}$, and $M$ is rigid when $k=0$ \cite[Proposition 2.4]{hpw2} and when $k=5$ \cite[Theorem 2]{crz}.
\end{proof}

\begin{proof}[Proof of Theorem \ref{thm_five}]
By \eqref{eq_quantization}, the possible values of $k$ are $0$, $2$, $3$ and $4$. If $k=0$, then $M$ is Einstein with $\Ric=\rho g=4\alpha g$, and it is isometric to the hemisphere $\Sph^5_+(\alpha)$ by \cite[Proposition 2.4]{hpw2}, where $\alpha=\lambda/(m+4)$ by \eqref{eq_constants}. If $k=2$ or $k=3$, the assertions are those of Theorem \ref{cor_alldim}, with $\alpha=\lambda/(m+2)$ and $\alpha=\lambda/(m+1)$, respectively. It remains to consider $k=4$. In this case, $R=4\lambda$ by \eqref{eq_quantization}, the universal cover of $M$ is compact \cite[Theorem 4.1]{hpw1}, and \cite[Theorem 2]{crz} shows that it is the product of a compact simply connected Einstein manifold $\widetilde N^4$ with $\Ric_{\widetilde N}=\lambda g_{\widetilde N}$ and an interval of length $\pi\sqrt{m/\lambda}$. In all four cases the universal cover is compact, and therefore $M$ is the quotient of it by a finite group.
\end{proof}

\section{Four-dimensional focal sets}\label{sec_four}

In this section, we prove Theorem \ref{thm_four}, so that $k=4$ and $F$ is a closed minimal submanifold of dimension four and codimension $q=n-4$. The Weyl tensor of a manifold of dimension two or three vanishes, which is no longer true in dimension four. For this reason, the identity of Theorem \ref{thm_general} now contains the pairing $W^F(\II,\II)$. This pairing has no sign a priori, and the term $2\alpha|\II|^2$ is not enough to control it, since the ratio $\alpha/\lambda=1/(m+q-1)$ is fixed.

\begin{proof}[Proof of Theorem \ref{thm_four}]
Taking $k=4$ in \eqref{eq_simons_weyl}, we have $\trg(\Ric_FT)+\omega_{\II}=-W^F(\II,\II)+2\inner{\mathring{\Ric}_F}{\mathring{T}}+\frac43\lambda t+\frac16t^2$. Recall that the integral of \eqref{eq_general_weyl} with $k=4$ vanishes by Theorem \ref{thm_general}, and that the integral of $I$ vanishes by Lemma \ref{lem_simons}. Thus the integral of $J-3I$ vanishes as well. Using $C_D=|D\II|^2-\frac12|A|^2$, $O_{\II}=\frac12\inner{K}{V}+C$, $|\mathring{T}|^2=|T|^2-\frac14t^2$ and
\begin{gather*}
	2\Big|V-\frac34K\Big|^2+\frac32\inner{K}{V}=2\Big|V-\frac38K\Big|^2+\frac{27}{16}C,\\
	6|\mathring{\Ric}_F|^2+2\inner{\mathring{\Ric}_F}{\mathring{T}}=6\Big|\mathring{\Ric}_F+\frac16\mathring{T}\Big|^2-\frac16|\mathring{T}|^2,
\end{gather*}
we arrive at
\begin{align}\label{eq_four_identity}
\begin{split}
	0=\int_F\Big(&|D\II|^2+\frac92|A|^2-3|d|^2+4\Big|L^s+\frac12\mathcal{P}\Big|^2+2\Big|V-\frac38K\Big|^2+2|Y|^2+\frac{1}{m^2\alpha^2}|\mathcal{F}|^2\\
	&+6\Big|\mathring{\Ric}_F+\frac16\mathring{T}\Big|^2+\Big(\frac23\lambda+2\alpha\Big)t-11W^F(\II,\II)+\frac{7}{16}C+\frac72c-\frac14t^2-\frac16|T|^2\Big)\dmu_F.
\end{split}
\end{align}
We will now estimate the terms of \eqref{eq_four_identity}. By \eqref{eq_d_bound} with $k=4$, we have $\frac92|A|^2-3|d|^2\geq0$. Moreover, the shape operators belong to the space of trace free symmetric endomorphisms of $T_pF$, which has dimension nine. Thus $G$ is a positive semidefinite Gram matrix of rank at most $\min\{q,9\}$ and trace $t$, which implies that $c\geq\frac19t^2$, and that $c\geq\frac12t^2$ if $q=2$. Now let $S$ be a trace free symmetric endomorphism of a four-dimensional space. Each eigenvalue of $S$ is minus the sum of the other three, and hence its square is at most three times the sum of their squares. It follows that $S^2\leq\frac34|S|^2\Id$. Consequently, $0\leq T\leq\frac34t\Id$ and $|T|^2\leq\frac34t\trg T=\frac34t^2$. Since $C\geq0$, we conclude that
\[
	\frac{7}{16}C+\frac72c-\frac14t^2-\frac16|T|^2\geq\Big(\frac{7}{18}-\frac14-\frac18\Big)t^2=\frac{1}{72}t^2,
\]
where $\frac{7}{18}$ is replaced by $\frac74$, and $\frac{1}{72}$ by $\frac{11}{8}$, if $q=2$, that is, if $n=6$. Discarding the remaining nonnegative terms in \eqref{eq_four_identity}, we get \eqref{eq_four_estimate}. Finally, suppose that \eqref{eq_four_pinching} holds. Then \eqref{eq_four_estimate} implies that $\int_F(|D\II|^2+\frac{1}{72}|\II|^4)\dmu_F\leq0$. Hence $\II=0$, and $M$ is rigid by Theorem \ref{thm_criterion}.
\end{proof}

\begin{remark}\label{rem_four}
Even among rigid manifolds, the maximum set can be any closed Einstein manifold $F^4$ with $\Ric_F=\lambda g_F$, because $F\times\Sph^q_+(\alpha)$ is a rigid manifold with maximum set $F$. In particular, constant scalar curvature does not imply that $W^F=0$, as the products with maximum set $\Sph^2\times\Sph^2$ or $\mathbb{CP}^2$ show. Nevertheless, the pairing $W^F(\II,\II)$ vanishes on every rigid manifold. It is natural to ask which additional hypotheses on the metric induced on $F$ give a sign to this pairing when $n=6$, in which case \eqref{eq_four_estimate} holds with the coefficient $\frac{11}{8}$.
\end{remark}

\appendix

\section{Two computations}\label{app_proofs}

In this appendix, we prove Lemmas \ref{lem_bochner_P} and \ref{lem_term_grad}, in which the parameter $m$ enters.

\begin{proof}[Proof of Lemma \ref{lem_bochner_P}]
Throughout this proof, we will write $\bar P=m\alpha g-P=\lambda g-\Ric$, $f_a=\nabla_af$ and $f_{ab}=\nabla^2_{a,b}f$. Since $\nabla f=-\frac{m}{u}\nabla u$, it follows from \eqref{fundeq} and \eqref{eq_P_identities} that
\[
	f_{ab}=\bar P_{ab}+\frac1mf_af_b,\qquad P_{ab}f_b=0,\qquad \bar P_{ab}f_b=m\alpha f_a,\qquad \nabla_aP_{ab}=\nabla_a\bar P_{ab}=0.
\]
The Ricci identity \eqref{eq_ricci_identity} gives $\nabla_af_{bc}-\nabla_bf_{ac}=R_{abcd}f_d$. Moreover, $\nabla_aP_{bc}=-\nabla_a\bar P_{bc}$ and
\[
	\nabla_a(f_bf_c)-\nabla_b(f_af_c)=f_bf_{ac}-f_af_{bc}=f_b\bar P_{ac}-f_a\bar P_{bc},
\]
which proves \eqref{eq_codazzi_P}. In index notation, this identity reads
\begin{align}\label{eq_codazzi_index}
	\nabla_aP_{bc}-\nabla_bP_{ac}=R_{abdc}f_d+\frac1m\big(f_b\bar P_{ac}-f_a\bar P_{bc}\big).
\end{align}
We will now compute $\Delta P$. Commuting derivatives by means of \eqref{eq_ricci_identity} and using $\nabla_aP_{ab}=0$, we get
\[
	\Delta P_{bc}=(\Ric P)_{bc}-\Rm(P)_{bc}+\nabla_a\big(\nabla_aP_{bc}-\nabla_bP_{ac}\big),
\]
so it remains to compute the divergence of the right-hand side of \eqref{eq_codazzi_index}. For its first term, the contracted second Bianchi identity and $P_{ab}f_b=0$ yield
\begin{align*}
	\nabla_a\big(R_{abdc}f_d\big)&=\nabla_{\nabla f}P_{bc}+P_{bd}f_{dc}+R_{abdc}f_{ad}\\
	&=\nabla_{\nabla f}P_{bc}+P_{bd}f_{dc}+m\alpha\Ric_{bc}-\Rm(P)_{bc}+\frac1mR_{abdc}f_af_d,
\end{align*}
and, for the second term,
\[
	\frac1m\nabla_a\big(f_b\bar P_{ac}-f_a\bar P_{bc}\big)=\frac1m\Big((\bar P^2)_{bc}+\alpha f_bf_c-(\Delta f)\bar P_{bc}+\nabla_{\nabla f}P_{bc}\Big).
\]
Let us now contract \eqref{eq_codazzi_index} with $f_a$. The result is
\[
	R_{abdc}f_af_d=\nabla_{\nabla f}P_{bc}+P_{bd}f_{dc}-\alpha f_bf_c+\frac{|\nabla f|^2}{m}\bar P_{bc}.
\]
Recalling Proposition \ref{prop_structure}(ii), we have $\Delta f=qm\alpha+\frac1m|\nabla f|^2$. It is also easy to see that $P_{bd}f_{dc}=m\alpha P_{bc}-(P^2)_{bc}$, and therefore $(\bar P^2)_{bc}+P_{bd}f_{dc}=m\alpha\bar P_{bc}$. Putting everything together, we see that the terms in $f_bf_c$ and in $|\nabla f|^2\bar P$ cancel, and it then follows that
\[
	\Delta P=\Ric P+m\alpha P-P^2+m\alpha\Ric-2\Rm(P)+(1-q)\alpha\bar P+\frac{m+2}{m}\nabla_{\nabla f}P.
\]
Now $(q-1)\alpha=\rho$ by \eqref{eq_constants}, and $\Ric=P+\rho g$. Thus the terms without derivatives add up to $(2m\alpha+2\rho)P-2\Rm(P)=2\lambda P-2\Rm(P)$. Since $\nabla f=-\frac{m}{u}\nabla u$, this is \eqref{eq_bochner_P}, as claimed. Finally, to prove \eqref{eq_bochner_global}, it suffices to take the inner product of \eqref{eq_bochner_P} with $P$ and to note that, since $|P|^2$ is constant, $\inner{\Delta P}{P}=-|\nabla P|^2$ and $\inner{\nabla_{\nabla u}P}{P}=0$.
\end{proof}

\begin{proof}[Proof of Lemma \ref{lem_term_grad}]
Let us arbitrarily fix $p\in F$. Recalling the Ricci identity \eqref{eq_ricci_identity}, a direct computation shows that
\begin{align*}
	\Delta\nabla_zP_{ab}&=\nabla_z\Delta P_{ab}+\Ric_{zx}\nabla_xP_{ab}\\
	&\quad+2\sum_{x,y}\big(R_{yzax}\nabla_yP_{xb}+R_{yzbx}\nabla_yP_{ax}\big)\\
	&\quad+\sum_x\big(\nabla_a\Ric_{zx}-\nabla_x\Ric_{za}\big)P_{xb}+\sum_x\big(\nabla_b\Ric_{zx}-\nabla_x\Ric_{zb}\big)P_{ax},
\end{align*}
where $\Ric_{zx}=\sum_yR_{yzyx}$. Note that the last line comes from the contracted second Bianchi identity $\sum_y\nabla_yR_{yzax}=\nabla_a\Ric_{zx}-\nabla_x\Ric_{za}$, and that it vanishes at $p$, since $\nabla\Ric=\nabla P$ is totally symmetric there (Lemma \ref{lem_first_order}(i)). We will now compute the first term on the right-hand side. With $f=-m\log u$, equation \eqref{eq_bochner_P} reads $\Delta P=2\lambda P-2\Rm(P)+\frac{m+2}{m}\nabla_{\nabla f}P$. Since $\nabla f=0$ and $\nabla^2f=m\alpha\,\Pi_N$ at $p$, we have $\nabla_z(\nabla_{\nabla f}P)=m\alpha(\Pi_N)_{zx}\nabla_xP$ there. On the other hand, $\Ric_{zx}=\lambda(\Pi_T)_{zx}+\rho(\Pi_N)_{zx}$ and, by \eqref{eq_constants}, $\rho+(m+2)\alpha=\lambda+2\alpha$. Therefore the sum of $\Ric_{zx}\nabla_xP_{ab}$ and $(m+2)\alpha(\Pi_N)_{zx}\nabla_xP_{ab}$ is $\lambda\nabla_zP_{ab}+2\alpha(\Pi_N)_{zx}\nabla_xP_{ab}$, whence it follows that
\begin{align*}
	\Delta\nabla_zP_{ab}&=3\lambda\nabla_zP_{ab}+2\alpha(\Pi_N)_{zx}\nabla_xP_{ab}-2(\nabla_z\Rm)(P)_{ab}-2\Rm(\nabla_zP)_{ab}\\
	&\quad+2\sum_{x,y}\big(R_{yzax}\nabla_yP_{xb}+R_{yzbx}\nabla_yP_{ax}\big).
\end{align*}
It remains to take the inner product of this identity with $\nabla_zP_{ab}$ and to sum over $z$. The first term gives $3\lambda|\nabla P|^2=9\lambda m^2\alpha^2t$, by Lemma \ref{lem_first_order}(i). Again by Lemma \ref{lem_first_order}(i), $\nabla_\beta P$ is $m\alpha S_\beta$ on $TF$ and vanishes elsewhere, so that the second term gives $2\alpha\sum_\beta|\nabla_\beta P|^2=2\alpha m^2\alpha^2t$. Since $P=m\alpha\,\Pi_T$, the third term is
\[
	-2m\alpha\sum(\nabla_z\Rm)_{ibid}(\nabla_zP)_{bd}=-2m^2\alpha^2\Omega,
\]
and the fourth one is $-2\Theta$. For the fifth term, we use the total symmetry of $\nabla P$ to write
\begin{align*}
	\sum R_{yzax}\nabla_yP_{xb}\nabla_zP_{ab}&=\sum_bR_{yzax}(\nabla P)_{yxb}(\nabla P)_{zab}\\
	&=-\sum_b\Rm\big((\nabla P)_{\cdot\cdot b},(\nabla P)_{\cdot\cdot b}\big)=-\Theta,
\end{align*}
where we used the skew symmetry of $R$ in its last two entries. Analogously, the last term gives $-\Theta$, and the result follows.
\end{proof}

\section*{AI Disclosure}

During the preparation of this manuscript, the author used publicly available large language models to discuss intermediate steps and to check computations involved in adapting the argument of Li, Li, Ou, Qu, Wu and Zhou \cite{lloqwz} to quasi-Einstein manifolds in Section \ref{sec_criterion}. No AI-generated text is present in this paper. The author independently verified all relevant mathematical arguments and takes full responsibility for the content of the manuscript.

\end{document}